\documentclass[a4paper]{article}

\usepackage{RR}
\usepackage{hyperref}

\usepackage[english]{babel}
\usepackage[latin1]{inputenc}
\usepackage[T1]{fontenc}
\usepackage{float}
\usepackage{amsmath}
\usepackage{amsfonts} 
\usepackage{amssymb} 
\usepackage{amsthm} 
\usepackage{stmaryrd} 
\usepackage{url} 
\usepackage{epsfig} 
\usepackage{natbib} 
\def\var{\mbox{var}} 
\def\corr{\textrm{corr}} 

\def\cov{\mbox{cov}}

\psfigdriver{dvips} 
 
\usepackage{mathrsfs} 
 
\usepackage{geometry} 
\newfloat{Table}{Hptb}{lot} 
 
\newtheorem{thrm}{Theorem} 
\newtheorem{prte}[thrm]{Proposition} 
\newtheorem{lemma}[thrm]{Lemma}

\RRdate{Novembre 2007}

\RRversion{3\thanks{Correction in Proposition 4.2}}
\RRdater{Mai 2008}
\RRNo{6354}
\RRauthor{Nicolas Verzelen
\thanks{Laboratoire de Math\'ematiques UMR 8628, Universit\'e Paris-Sud, 91405 Osay}
\thanks{INRIA Futurs, Projet SELECT, Universit\'e Paris-Sud, 91405 Osay}%
  \and Fanny Villers \thanks{INRA, Math\'ematiques et Informatique Appliqu\'ees MIA, 78352 Jouy-en-Josas}%
}
\authorhead{Verzelen \& Villers}
\RRtitle{Tests d'hypoth\`eses lin\'eaires pour des mod\`eles de r\'egression gaussienne en grande dimension}
\RRetitle{Goodness-of-fit Tests for high-dimensional Gaussian linear models}
\titlehead{Goodness-of-fit Tests for high-dimensional Gaussian linear models}
\RRresume{On consid\'ere $(Y,(X_i)_{i\in\mathcal{I}})$ un vecteur gaussien de moyenne nulle et $V$ un sous-ensemble de $\mathcal{I}$. Supposons que l'on observe $n$ r\'eplications ind\'ependantes du vecteur $(Y,X)$. Dans ce rapport, nous proposons une proc\'edure pour tester l'hypoth\`ese $Y$ est ind\'ependant de
$(X_i)_{i\in \mathcal{I}\backslash V}$ conditionellement \`a $(X_i)_{i\in V}$. Ce test ne n\'ecessite aucune connaissance sur la matrice de covariance de $X$ ou la variance de $Y$ et peut s'appliquer dans un contexte de grande dimension. De plus, on peut facilement en d\'eduire un test de voisinage de mod\`ele graphique gaussien. Nous calculons la puissance du test d'un point de vue non asymptotique et prouvons qu'il atteint la vitesse optimale de s\'eparation (\`a un facteur $\log n$ pr\'es) pour diff\'erentes classes d'alternatives. Nous en d\'eduisons ainsi des vitesse minimax de s\'eparation d'hypoth\`eses pour ce mod\`ele de r\'egression. Enfin, nous \'evaluons la performance de notre proc\'edure sur des donn\'ees simul\'ees.
}
\RRabstract{Let $(Y,(X_i)_{1\leq i\leq p})$ be a real zero mean Gaussian vector and $V$ be a
subset of $\{1,\ldots, p \}$. Suppose we are given $n$ i.i.d. replications of this vector. We propose a new test for testing that $Y$ is independent of
$(X_i)_{i\in \{1,\ldots p\}\setminus V}$ conditionally to $(X_i)_{i\in V}$ against the
general alternative that it is not. This procedure does not depend on any prior information on
the covariance of $X$ or the variance of $Y$ and applies in a high-dimensional
setting. It straightforwardly extends to test the neighbourhood of a Gaussian graphical model. 
The procedure is based on a model of Gaussian regression with random Gaussian covariates. We give non asymptotic properties of the test and
we prove that it is rate optimal (up to a possible $\log(n)$ factor) over various classes of alternatives
under some additional assumptions. Besides, it allows us to
derive non asymptotic minimax rates of testing in this random design setting. Finally, we
carry out a simulation study in order to evaluate the performance of our procedure.    }
\RRmotcle{R\'egression lin\'eaire, mod\`eles graphiques gaussiens, test multiple, test adaptatif, test minimax, vitesse de s\'eparation minimax, ellipso\"ide, test d'ad\'equation}
\RRkeyword{Linear regression, Gaussian graphical models,
multiple testing, adaptive testing, minimax hypothesis testing, minimax
separation rate, ellipsoid, goodness-of-fit}
\RRprojet{Select}
\RRtheme{\THCog} 
\URFuturs 
\begin{document}
\makeRR   

\section{Introduction}  
We consider the following regression model 
\begin{eqnarray}\label{model} 
Y = \sum_{i=1}^p\theta_iX_i +\epsilon 
\end{eqnarray} 
where $\theta$ is an unknown vector of $\mathbb{R}^{p}$. In the sequel, we note $\mathcal{I}:=\left\{1,\ldots,p\right\}$. The vector 
$X:=(X_i)_{1\leq i\leq p}$ follows a real zero mean Gaussian distribution with non singular covariance 
matrix $\Sigma$ and $\epsilon$ is a real zero mean Gaussian random variable 
independent of $X$. Straightforwardly, the variance of $\epsilon$ 
corresponds to the conditional variance of $Y$ given $X$, $\var(Y|X)$.\\

The variable selection problem for this model in a high-dimensional setting has recently attracted a lot of attention. A large number of papers are now devoted to the design of new algorithms and estimators which are computationally feasible and are proven to converge; see for instance the works of Meinshausen and B\"uhlmann \cite{meinshausen06}, Cand\`es and Tao \cite{candes07}, Zhao and Yu \cite{zhao06}, Zou and Hastie \cite{zou05}, B\"uhlmann and Kalisch \cite{buhlmann07}, or Zhao and Huang \cite{zhang08}. A common drawback of the previously mentioned estimation procedures is that they require  restrictive conditions on the covariance matrix $\Sigma$ in order to behave well. Our issue is the natural testing counterpart of this variable selection problem: we aim at defining a computationally feasible testing procedure that achieves an optimal rate for any covariance matrix $\Sigma$.
 
\subsection{Presentation of the main results}

We are given $n$ i.i.d. replications of the vector $(Y,X)$. Let us respectively note ${\bf Y}$ and ${\bf 
  X}_i$ the vectors of the $n$ observations of $Y$ and 
$X_i$ for any $i\in\mathcal{I}$. Let $V$ be a subset of $\mathcal{I}$, then $X_V$ refers to the set $\left\{X_i,i\in 
V\right\}$ and $\theta_V$ stands for the sequence $(\theta_i)_{i\in V}$. We first propose a collection of testing procedures $T_{\alpha}$ of the null hypothesis 
``$\theta_{\mathcal{I} \setminus  V} 
= 0$'' against the general alternative ``$\theta_{\mathcal{I} \setminus  V} 
\neq 0$''. These procedures are based on the ideas of Baraud \emph{et al.} \cite{baraud2003} in a random design. Their definition are very flexible as they require no prior knowledge of the covariance of $X$, the 
variance of $\epsilon$, nor the variance of $Y$. Note that the property ``$\theta_{\mathcal{I} \setminus  V} 
= 0$'' is equivalent to ``$Y$ is independent of
  $X_{\mathcal{I}\backslash V}$ conditionally to $X_V$''. Hence, it also permits to test conditional independences and applies for testing the graph of Gaussian graphical model (see below). Contrary to most approaches in this setting (e.g. Drton and Pearlman \cite{drton07}), we are able to consider the difficult case of tests in a high-dimensional setting: the number of covariates $p$ is possibly much larger than the number of observations $n$. Such situations arise in many statistical applications like in genomics or biomedical imaging. To our knowledge, the only testing procedures (e.g. \cite{schafer05}) that could handle high-dimensional alternatives lack of theoretical justifications. In this paper, we exhibit some tests $T_{\alpha}$ that are both computationally amenable and optimal in the minimax sense.

From a theoretical perspective, we are able to control the Family Wise Error Rate (FWER) of our testing procedures $T_{\alpha}$. Besides, we derive a general non asymptotic upper bound for their power. Contrary to the various rates of convergence obtained in the estimation setting (e.g. \cite{meinshausen06} or \cite{candes07}), our upper bound holds for any covariance matrix $\Sigma$. Then, we derive from it non-asymptotic minimax rates of testing in the Gaussian random design framework. If the minimax rates are known for a long time in the fixed design Gaussian regression framework (e.g.\cite{baraudminimax}), they were unknown in our setting. For instance, if at most $k$ components of $\theta$ are non-zero and if $k$ is much smaller than $p$, we prove that the minimax rates of testing is of order $\frac{k\log(p)}{n}$ when the covariates $X_i$ are independent. If the covariates are dependent, we derive faster minimax rates. To our knowledge, these are the first results for testing  or estimation issues that illustrate minimax rates for dependent covariates. Afterwards, we show analogous results when $k$ is large, or when the vector $\theta$ belongs to some ellipsoid or some collection of ellipsoids. For any of these alternatives, we exhibit some procedure $T_{\alpha}$ that achieves the optimal rate (at a possible $\log(n)$ factor). Finally, we illustrate the performance of the procedure on simulated examples.

\subsection{Application to Gaussian Graphical Models (GGM)} 
 
Our work was originally motivated by the following question: let 
$(Z_j)_{j\in\mathcal{J}}$ be a random vector which follows a zero 
mean Gaussian distribution whose covariance matrix $\Sigma'$ is non singular. We observe $n$ 
i.i.d. replications of this vector $Z$ and we are given a graph 
$\mathcal{G}=(\Gamma,E)$ where  $\Gamma=\{1,\ldots |\mathcal{J}|\}$ and $E$ is a set of edges in $\Gamma\times \Gamma$. How can we test that $Z$ is an undirected Gaussian graphical model (GGM) with 
respect to the graph $\mathcal{G}$?

The random vector $Z$ is a GGM with respect to the graph $\mathcal{G}=(\Gamma,E)$ if for any couple $(i,j)$ which is not contained in the edge set $E$, $Z_i$ and $Z_j$ are independent, given the remaining variables. See Lauritzen \cite{lauritzen96} for definitions 
and main properties of GGM.  Interest in these models has grown as they allow the description of dependence structure in high-dimensional data. As such, they are widely used in spatial statistics \citep{cressie,rue} or probabilistic expert systems \citep{cowell}. More recently, they have been applied to the analysis of microarray data. The challenge is to infer the network regulating the expression of the genes using only a small sample of data, see for instance Sch\"afer and Strimmer \cite{schafer05}, Kishino and Waddell \cite{kishino00}, or Wille \emph{et al.} \cite{wille04}. This issue has motivated the research for new estimation procedures to handle GGM in a high-dimensional setting. 

It is beyond the scope of this paper to give an exhaustive review of these. Many of these graph estimation methods are based on multiple testing procedures, see for instance Sch\"afer and Strimmer \cite{schafer05} or Wille and B\"uhlmann \cite{wille06}. Other methods are based on variable selection for high-dimensional data we previously mentioned. For instance, Meinshausen and B\"uhlmann \cite{meinshausen06} proposed a computationally feasible model selection algorithm using Lasso penalization. Huang \emph{et al.} \cite{huang06} and Yuan and Lin \cite{yuan07} extend this method to infer directly the inverse covariance matrix $\Sigma'^{-1}$ by minimizing the log-likehood penalized by the $l^1$ norm. 

While the issue of graph and covariance estimation is extensively studied, few theoretical results are proved for the problem of hypothesis testing of GGM in a high-dimensional setting. We believe that this issue is significant for two reasons: first, when considering a gene regulation network, the biologists often have a previous knowledge of the graph and may want to test if the microarray data match with their model. Second, when applying an estimation method in a high-dimensional setting, it could be useful to test the estimated graph as some of these methods reveal too conservative.

Admittedly, some of the previously mentioned estimation methods are based on multiple testing. However, as they are constructed for an estimation purpose, most of them do not take into account some previous knowledge about the graph. This  
is for instance the case for the approaches of Drton and Perlman \cite{drton07} and Sch\"afer and Strimmer \cite{schafer05}. Some of the other existing procedures cannot be applied in a high-dimensional setting ($ |\mathcal{J}|\geq n$). Finally, most of them lack theoretical justification in a non asymptotic way.

In a subsequent paper \citep{villers_ggm} we define   
a test of graph based on the present work. It benefits the ability of handling high 
dimensional GGM and has minimax properties. Besides we show numerical 
evidence of its efficiency; see \cite{villers_ggm} for more details. In this article, we shall only present the idea underlying our approach. 
 
For any $j\in \mathcal{J}$, we note $N(j)$ the set of neighbours of $j$ in the 
graph $\mathcal{G}$. Testing 
that $Z$ is a GGM with respect to $\mathcal{G}$ is equivalent to testing  
that the random variable $Z_j$ conditionally to 
$(Z_l)_{l\in N(j)}$ is independent of 
$(Z_l)_{l\in\mathcal{J} \setminus  (N(j)\cup \{j\})}$ for any $j\in 
\mathcal{J}$. As $Z$ follows a Gaussian distribution, the distribution of $Z_j$ conditionally to the other 
variables decomposes as follows: 
\begin{eqnarray*} 
Z_j  = \sum_{k\in \mathcal{J} \setminus \{j\}} \theta_k Z_k +\epsilon_j, 
\end{eqnarray*} 
where $\epsilon_j$ is normal and independent of $(Z_k)_{k\in \mathcal{J} \setminus
 \{j\}}$. Then, the statement of conditional independency is equivalent to 
 $\theta_{\mathcal{J} \setminus  \{j\}\cup N(j)}=0$.  
 This approach based on conditional 
 regression is also used for estimation by Meinshausen and B\"uhlmann \cite{meinshausen06}.\\

\subsection{Organization of the Paper} 
 
In Section \ref{section_description}, we present the approach of our procedure and connect it with the fixed design framework. Besides, we define the notion of minimax rates of testing in this setting and gather the main notations.
We define the testing procedures $T_{\alpha}$ in Section \ref{section_definition} and we non asymptotically characterise  
the set of vectors $\theta$ over which the test $T_{\alpha}$ is powerful. In Section \ref{detection} and \ref{ellipsoides}, we apply our procedure to 
define tests and study their optimality for two different classes of 
alternatives. More precisely, in Section \ref{detection} we test $\theta=0$ against 
the class of $\theta$ whose components equal $0$, except at most 
$k$ of them ($k$ is supposed small). We define a test which under mild conditions achieves the minimax 
rate of testing. When the covariates are independent, it is interesting to 
note that the minimax rates exhibits the same ranges in our statistical model 
(\ref{model}) and in fixed design regression model 
(\ref{model_fixed_design}). In Section \ref{ellipsoides}, we define two 
procedures that achieve the simultaneous minimax rates of testing over large 
classes of ellipsoids (to sometimes the price of a $\log(p)$ factor). Besides, we show that 
the problem of adaptation over classes of ellipsoids is impossible without a 
loss in efficiency. This was previously pointed out in \cite{spokoiny96} in 
fixed design regression framework. The 
simulation studies are presented in Section \ref{simulations1}. Finally, Sections 
\ref{proof1}, \ref{proof2} and Appendix contain the proofs.

\section{Description of the approach}\label{section_description}

\subsection{Connection with tests in fixed design regression} 
 
Our work is directly inspired by the testing procedure of Baraud \emph{et al.} \cite{baraud2003} in 
fixed design regression framework. Contrary to model (\ref{model}), the problem of hypothesis testing in fixed 
design regression has 
been extensively studied. This is why we 
will use the results in this framework as a benchmark for the theoretical 
bounds in our model (\ref{model}). Let us define this second regression model: 
\begin{eqnarray}\label{model_fixed_design} 
Y_i = f_i +\sigma \epsilon_i,\ \  i\in\{1,\ldots,N\}, 
\end{eqnarray} 
where $f$ is an unknown vector of $\mathbb{R}^N$, $\sigma$ some unknown 
positive number, and the $\epsilon_i$'s a sequence of i.i.d. standard Gaussian 
random variables. The problem at hand is testing that $f$ belongs to a linear subspace of 
$\mathbb{R}^N$ against the alternative that it does not. We refer to \cite{baraud2003} for a short review of non parametric 
tests in this framework. Besides, we are interested in the performance of the 
procedures from a minimax perspective. To our knowledge, there has been no 
results in model (\ref{model}). However, there are numerous papers 
on this issue in the fixed design regression model. First, we refer to the 
seminal work of Ingster \citep{ingster93a,ingster93b,ingster93c} who gives asymptotic minimax rates over non 
parametric alternatives. Our work is closely related to the results of 
Baraud \cite{baraudminimax} where he gives non asymptotic minimax rates of testing over 
ellipsoids or sparse signals. Throughout the paper, we highlight the link 
between the minimax rates in fixed and in random design. \\ 
 
\subsection{Principle of our testing procedure}\label{presentation} 
 
Let us briefly describe the idea underlying our testing procedure. A formal definition will follow in Section \ref{section_definition_Procedure}.Let $m$ be a 
subset of $\mathcal{I} \setminus  V$.  We respectively define $S_V$ and $S_{V\cup m}$  
as the linear subspaces of $\mathbb{R}^{p}$ such that 
$\theta_{\mathcal{I} \setminus  V} = 0$, respectively $\theta_{\mathcal{I} \setminus  
  (V\cup m)} = 0$. We note $d$ and $D_m$ for the cardinalities of $V$ 
and $m$ and $N_m$ refers to $N_m=n-d-D_m$. 
If $N_m>0$, we define the Fisher statistic $\phi_m$ by 
\begin{eqnarray}\label{definition_phi} 
\phi_m({\bf Y},{\bf X}) := \frac{N_m\|\Pi_{V\cup m}{\bf Y} - 
  \Pi_{V}{\bf{Y}}\|_n^2}{D_m\|{\bf Y} - \Pi_{V\cup m}{\bf Y}\|^2_n}\ , 
\end{eqnarray} 
where $ \Pi_{V}$ refers to the orthogonal projection onto the space generated by 
the vectors $({\bf X}_i)_{i\in V}$ and $\|.\|_n$ is the canonical norm in $\mathbb{R}^n$.
 We define the test statistic $\phi_{m,\alpha}({\bf Y},{\bf X})$ as 
\begin{eqnarray}\label{definition_phi_test} 
\phi_{m,\alpha}({\bf Y},{\bf X}) = \phi_m({\bf 
  Y},{\bf X}) - \bar{F}^{-1}_{D_m,N_m}(\alpha)\ ,
\end{eqnarray} 
 where $\bar{F}_{D_m,N_m}(u)$ denotes the probability for a Fisher variable with $D$ and $N$ degrees of freedom to be larger than $u$.
Let us consider a finite collection $\mathcal{M}$ of non empty 
subsets of $\mathcal{I} \setminus  V$ such that for each $m\in\mathcal{M}$, $N_m>0$. Our 
testing procedure consists of doing a Fisher test for each 
$m\in\mathcal{M}$. We define $\{\alpha_m,m\in\mathcal{M}\}$ a suitable 
collection of numbers in $]0,1[$ (which possibly depends on ${\bf X}$). For 
    each $m\in\mathcal{M}$, we do the Fisher test $\phi_m$ of level $\alpha_m$ of: 
    $$H_0:\theta\in S_V \text{\ \ against the alternative \ 
      \ } H_{1,m}: \theta\in S_{V\cup m}\setminus S_V$$ 
    and we decide to reject the null hypothesis if one of those Fisher tests does.\\ 
 
The main advantage of our procedure is that it is very flexible in the choices 
of the model $m\in\mathcal{M}$ and in the choices of the weights 
$\{\alpha_m\}$. Consequently, if we choose a suitable collection 
$\mathcal{M}$, the test is powerful over a large class of 
alternatives as shown in Sections \ref{section_puissance}, \ref{detection},
and \ref{ellipsoides}. \\ 

Finally, let us mention that our procedure easily extends to the
case where the expectation of the random vector $(Y,X)$ is unknown. Let ${\bf \overline{X}}$ and ${\bf \overline{Y}}$ denote the projections of  ${\bf X}$ and ${\bf Y}$ onto the unit vector ${\bf 1}$. Then, one only has to
apply the procedure to $({\bf Y-\overline{Y}}, {\bf X-\overline{X}})$ and to
replace $d$ by $d+1$. The properties of the test remain unchanged and one can adapt all the proofs to the price of more technicalities.

\subsection{Minimax rates of testing} 
 
In order to examine the quality of our tests, we will compare their 
performance with the minimax rates of testing. That is why we now define precisely what we mean by the $(\alpha, \delta)$-minimax rate of 
testing over a set $\Theta$. We endow $\mathbb{R}^{p}$ with the Euclidean norm
\begin{eqnarray}\label{definition_produit}
 \|\theta\|^2 := \theta^t\Sigma \theta= \var\left(\sum_{i=1}^p\theta_i X_i\right)\ .	
\end{eqnarray}

As $\epsilon$ and $X$ are independent, we derive from the definition of $\|.\|^2$ that 
$\var(Y) = \|\theta\|^2 + \var(Y|X)$. Let us remark that $\var(Y|X)$ does not depend on $X$. If we have $\|\theta\|$ vary, either
the quantity $\var(Y)$ or $\var(Y|X)$ has to vary. In the sequel, we suppose that $\var(Y)$ is
fixed. We briefly justify this choice in Section \ref{section_independent}. Consequently, if $\|\theta\|^2$ is increasing, then $\var(Y|X)$ has to
decrease so that the sum remains constant. Let $\alpha$ be a number in $]0;1[$ and
let $\delta$ be a number in $]0;1-\alpha[$ (typically small).
 For a given vector $\theta$, matrix $\Sigma$ and $\var(Y)$, we denote 
$\mathbb{P}_{\theta}$ the joint distribution of $({\bf Y},{\bf X})$. For the 
    sake of simplicity, we do not emphasize the dependence of 
    $\mathbb{P}_{\theta}$ on $\var(Y)$ or  $\Sigma$. Let $\psi_{\alpha}$ be a test of level $\alpha$ of the hypothesis "$\theta=0$" against the hypothesis "$\theta\in\Theta\setminus 0$".
In our framework, it is natural to measure the 
performance of $\psi_{\alpha}$ using the quantity 
$\rho\left(\psi_{\alpha},\Theta,\delta,\var(Y),\Sigma\right)$ defined by: 
\begin{eqnarray*} 
\rho\left(\psi_{\alpha},\Theta,\delta,\var(Y),\Sigma\right) & := & 
\inf\left\{\rho>0,\ \inf\left\{\mathbb{P}_{\theta}(\psi_{\alpha}=1),\ \theta\in\Theta \text{ and } \frac{\|\theta\|^2}{\var(Y)-\|\theta\|^2}\geq
  \rho^2\right\}\geq 1-\delta \right\}, 
\end{eqnarray*}  
where the quantity
\begin{eqnarray}\label{ratio} 
r_{s/n}(\theta):=\frac{\|\theta\|^2}{\var(Y)-\|\theta\|^2} 
\end{eqnarray} appears naturally as 
it corresponds to the ratio $\|\theta\|^2/\var(Y|X)$ which is  the 
quantity of information brought by $X$ (i.e. the signal) over the  conditional variance of 
$Y$ (i.e. the noise). 
We aim at describing the quantity  
\begin{eqnarray} 
\inf_{\psi_{\alpha}}\rho\left(\psi_{\alpha},\Theta,\delta,\var(Y),\Sigma\right) := \rho\left(\Theta,\alpha,\delta,\var(Y),\Sigma\right), 
\end{eqnarray} 
where the infimum is taken over all the level-$\alpha$ tests 
$\psi_{\alpha}$. We call this quantity the $(\alpha,\delta)$-minimax rate of 
testing over $\Theta$. 
 
A dual notion of this $\rho$ function is the function $\beta_{\Sigma}$. For any 
$\Theta\subset \mathbb{R}^p$ and $\alpha\in ]0,1[$, we denote 
    $\beta_\Sigma(\Theta)$ the quantity 
\begin{eqnarray*} 
\beta_{\Sigma}\left(\Theta\right) := \inf_{\psi_{\alpha}}\sup_{\theta\in \Theta}\mathbb{P}_{\theta}\left[\psi_{\alpha}=0\right], 
\end{eqnarray*} 
where the infimum is taken over all level-$\alpha$ tests $\psi_{\alpha}$ and 
where we recall that $\Sigma$ refers to the covariance matrix of $X$.\\

\subsection{Notations}
 
Let recall the main notations that we shall use throughout the 
paper. In the sequel, $n$ stands for the number of independent observations, $p$ is the number of covariates. Besides, $X_V$ stands for the collection $(X_{i})_{i\in V}$ of the covariates that correspond to the null hypothesis and $d$ is the cardinality of the set $V$. The models $m$ are subsets of $\mathcal{I}\subset V$ and we note $D_m$ their cardinality. $T_{\alpha}$ stands for our testing procedure of level $\alpha$. The statistics $\phi_m$ and the test $\phi_{m,\alpha}$ are respectively defined in (\ref{definition_phi}) and (\ref{definition_phi_test}). Finally, the norm $\|.\|$ is introduced in \ref{definition_produit}.

For $x,y\in \mathbb{R}$, we set 
$$x\wedge y := \inf\{x,y\},\, x\vee y := \sup\{x,y\}.$$ 
For any $u\in \mathbb{R}$, $\bar{F}_{D,N}(u)$ denotes the probability for a 
Fisher variable with $D$ and $N$ degrees of freedom to be larger than $u$. In the sequel, $L$, $L_1$, $L_2$,$\ldots$ denote constants that may vary from line to line. The notation $L(.)$ specifies the dependency on some quantities. For the sake of simplicity, we only give the orders of magnitude in the results and we refer to the proofs for explicit constants.

\section{The Testing procedure}\label{section_definition}

\subsection{Description of the procedure} \label{section_definition_Procedure}
 
Let us first fix some level $\alpha\in ]0,1[$. Throughout this paper, we suppose 
    that $n\geq d+2$. Let us consider a finite collection $\mathcal{M}$ of non 
    empty subsets of $\mathcal{I}\setminus V$ such that for all $m\in \mathcal{M}$, 
    $1\leq D_m\leq n-d-1$. We introduce the following test of level $\alpha$. We reject $H_0$: ``$\theta\in S_V$'' when the statistic 
\begin{eqnarray} 
T_{\alpha} := \sup_{m\in\mathcal{M}}\left\{\phi_m({\bf Y},{\bf X}) - \bar{F}^{-1}_{D_m,N_m}(\alpha_m({\bf X}))\right\}\label{definition_test} 
\end{eqnarray} 
is positive, where the collection of weights $\left\{\alpha_m({\bf X}),m\in\mathcal{M}\right\}$ is chosen according to one of the two following procedures:\vspace{0.5cm}~\\ 
$P_1$ : The $\alpha_m$ 's do not depend on ${\bf X}$ and satisfy the equality :  
\begin{equation} 
\sum_{m \in \mathcal{M} } \alpha_m = \alpha \ .
\end{equation}~\\ 
$P_2$ :  For all $m \in \mathcal{M}$, $\alpha_m ({\bf X})= q_{{\bf X}, \alpha}$, the $\alpha$-quantile of the 
distribution of the random variable 
\begin{eqnarray}  
\inf_{m \in \mathcal{M}} \bar{F}_{D_m,N_m}\left( \frac{ \| 
  \Pi_{V\cup m}(\boldsymbol{\epsilon}) -  \Pi_{V}(\boldsymbol{\epsilon})  
  \|_n^2 /D_m }{ \|\boldsymbol{\epsilon} - 
  \Pi_{V\cup m}(\boldsymbol{\epsilon}) \|_n^2 /N_m}\right)\label{methode_quantile} 
\end{eqnarray} 
conditionally to ${\bf X}$.

Note that it is easy to compute the quantity $q_{{\bf X},\alpha}$. Let $Z$ be 
a standard Gaussian random vector of size $n$ 
independent of ${\bf X}$. As ${\boldsymbol{\epsilon}}$ is independent of 
${\bf X}$, the distribution of (\ref{methode_quantile}) conditionally to ${\bf 
  X}$ is the same as the 
distribution of  
$$\inf_{m \in \mathcal{M}} \bar{F}_{D_m,N_m} \left(\frac{ \| 
  \Pi_{V\cup m}(Z) -  \Pi_{V}(Z)  
  \|^2 /D_m }{ \|Z - 
  \Pi_{V\cup m}(Z) \|^2 /N_m}\right)$$ 
conditionally to ${\bf X}$. Hence, we can easily work out its quantile using 
  Monte-Carlo method. \vspace{0.5cm}
 
Clearly, the computational complexity of the procedure is linear with respect to the size of the collection of models $\mathcal{M}$ even when using Procedure $P_2$. Consequently, when we apply our procedure to high-dimensional data as in Section \ref{simulations1} or in \cite{villers_ggm}, we favour collections $\mathcal{M}$ whose size is linear with respect to the number of covariates $p$.  

\subsection{Comparison of Procedures $P_1$ and $P_2$}\label{compa_puissance} 

We respectively refer to $T^1_{\alpha}$ and $T^2_{\alpha}$ for the tests (\ref{definition_test}) 
associated with Procedure $P_1$ and $P_2$. First, we are able to control the behavior of the test under the null hypothesis. 

\begin{prte}\label{niveau}
The test $T^1_{\alpha}$ corresponds to a Bonferroni procedure and therefore satisfies 
\begin{eqnarray*}
 \mathbb{P}_{\theta}(T_{\alpha} > 0) 
 \leq \sum_{m \in \mathcal{M}}\alpha_m\leq \alpha,
\end{eqnarray*}
whereas the test $T^2_{\alpha}$ has the property to be of size exactly $\alpha$:
\begin{eqnarray*}
 \mathbb{P}_{\theta}(T_{\alpha} > 0) = \alpha.
\end{eqnarray*}
\end{prte}
The proof is given in Appendix. Besides, the test $T^2_{\alpha}$ is more powerful
 than the corresponding test $T^1_{\alpha}$ defined with weights $\alpha_m= \alpha/|\mathcal{M}|$.

\begin{prte}\label{puissance_comparaison}
For any parameter $\theta$ that does not belong to $S_V$, the procedure $T^1_{\alpha}$ with weights $\alpha_m= \alpha/|\mathcal{M}|$ and the procedure $T^2_{\alpha}$ satisfy
\begin{eqnarray}\label{comparaison_puissance} 
\mathbb{P}_{\theta}\left(\left. T^2_{\alpha}({\bf X},{\bf Y}) >0\right|{\bf X}\right)\geq 
\mathbb{P}_{\theta}\left(T^1_{\alpha}({\bf X},{\bf Y})>0\left.\right|{\bf 
  X}\right)\,\, \, \, {\bf X}\ a.s.\ . 
\end{eqnarray}  
\end{prte}
Again, the proof is given in Appendix.
On the one hand, the choice of Procedure $P_1$ allows to avoid the computation of 
the quantile $q_{{\bf X},\alpha}$ and possibly permits to give a 
Bayesian flavor to the choice of the weights. On the other hand, Procedure $P_2$ is more powerful 
than the corresponding test with Procedure $P_1$. We will illustrate these considerations
in Section \nolinebreak\ref{simulations1}. In  sections \ref{section_puissance}, \ref{detection}, and \ref{ellipsoides} we study the power 
and rates of testing of $T_{\alpha}$ with Procedure $P_1$.

\subsection{Power of the Test}\label{section_puissance} 
 
We aim at describing a set of vectors $\theta$ in $\mathbb{R}^{p}$ over which the test defined 
in Section \ref{section_definition} with Procedure $P_1$ is powerful. Since Procedure $P_2$ is more powerful than Procedure $P_1$ with $\alpha_m=\alpha/|\mathcal{M}|$, the test with Procedure $P_2$ will also be powerful on this set of \nolinebreak$\theta$.

Let  $\alpha$ and $\delta$ be two numbers in $]0,1[$, and let $\{\alpha_m, m\in 
    \mathcal{M}\}$ be weights such that $\sum_{m\in\mathcal{M}}\alpha_m\nolinebreak\leq \nolinebreak\alpha$. Let define Hypothesis $(H_{\mathcal{M}})$ as follows:~\\~\\
$(H_{\mathcal{M}})$ \ \ \ \ \ \ \ \ \ \ \ \ \    For all $m\in \mathcal{M}$, $\alpha_m\geq \exp(-N_m/10)$ and $\delta\geq \exp 2(-N_m/21)$.~\\ 

For typical choices of the collections $\mathcal{M}$ and 
$\{\alpha_m,m\in\mathcal{M}\}$, these conditions are fulfilled as  
discussed in Sections \ref{detection} and \ref{ellipsoides}. Let us now turn to the main result. 
 
\begin{thrm}\label{thrm_puissance} 
Let $T_{\alpha}$ be the test procedure defined by (\ref{definition_test}). We assume 
that $n>d+2$ and that Assumption $(H_{\mathcal{M}})$ holds. Then, $\mathbb{P}_{\theta}(T_\alpha>0)\geq 1-\delta$ for all 
$\theta$ belonging to the set 
$$\mathcal{F}_{\mathcal{M}}(\delta):=\left\{\theta\in 
\mathbb{R}^{p},\exists m\in  \mathcal{M}: 
\frac{\emph{var}(Y|X_V)-\emph{var}(Y|X_{V\cup m})}{\emph{var}(Y|X_{V\cup m})} \geq \Delta(m)\right\},$$ 
where 
\begin{eqnarray} 
\Delta(m) := \frac{L_1\sqrt{D_m\log\left(\frac{2}{\alpha_m \delta}\right)}\left(1+\sqrt{\frac{D_m}{N_m}}\right)+L_2 \left(1+2\frac{D_m}{N_m} \right)\log\left(\frac{2}{\alpha_m \delta}\right)}{n-d}\label{inegalite_puissance_hm}.
\end{eqnarray} 

\end{thrm} 
 
This result is similar to Theorem 1 in \cite{baraud2003} in fixed design 
regression framework and the same comment 
also holds: the test $T_{\alpha}$ under procedure $P_1$ has a power comparable 
to the best of the tests among the family $\left\{\phi_{m,\alpha}, 
m\in\mathcal{M}\right\}$. Indeed, let us assume for instance that $V=\{0\}$ and that the $\alpha_m$ are chosen to be equal to 
$\alpha/|\mathcal{M}|$. The test $T_{\alpha}$ defined by 
 (\ref{definition_test}) is equivalent to doing several tests of $\theta=0$ against 
$\theta\in S_m$ at level $\alpha_m$ for $m\in\mathcal{M}$ and it rejects the null hypothesis if one of 
those tests does. From Theorem \ref{thrm_puissance}, we know that under the hypothesis $H_{\mathcal{M}}$ this test has 
a power greater than $1-\delta$ over the set of vectors $\theta$ belonging to 
$\bigcup_{m\in\mathcal{M}}\mathcal{F}'_m(\delta,\alpha_m)$ where $\mathcal{F}'_m(\delta,\alpha_m)$ is the set of vectors $\theta\in \mathbb{R}^{p}$ such that
\begin{eqnarray}\label{puissance_multitest} 
\frac{\var(Y)- \var(Y|X_m)}{\var(Y|X_m)}\geq \frac{L(D_m,N_m)}{n}\left(\sqrt{D_m\log\left(\frac{2}{\alpha_m\delta}\right)}+\log\left(\frac{2}{\alpha_m\delta}\right)\right)\ . 
\end{eqnarray} 
Besides, $L(D_m,N_m)$ behaves like a constant if the ratio $D_m/N_m$ is 
bounded. Let us compare this result with the set of $\theta$ over which the 
Fisher test $\phi_{m,\alpha}$ at level $\alpha$ has a power greater than 
$1-\delta$. Applying Theorem \ref{thrm_puissance}, we know that it contains $\mathcal{F}'_m(\delta,\alpha)$. 
Moreover, the following Proposition shows that it is not much larger than $\mathcal{F}'_m(\delta,\alpha)$: 
 
\begin{prte}\label{minoration_puissance_1test} 
Let $\delta\in ]0,1-\alpha[$. 
If  
$$\frac{\emph{var}(Y)-\emph{var}(Y|X_m)}{\emph{var}(Y|X_m)}\leq L(\alpha,\delta)\frac{\sqrt{D_m}}{n},$$ 
then $\mathbb{P}_{\theta}\left(\phi_{m,\alpha} >0\right)\leq 1-\delta$  
\end{prte} 
The proof is postponed to Section \ref{proof2} and is based on a lower bound 
of the minimax rate of testing.

$\mathcal{F}'_m(\delta,\alpha)$ and $\mathcal{F}'_m(\delta,\alpha_m)$ defined 
in  (\ref{puissance_multitest}) differ from the fact that $\log(1/\alpha)$ is replaced by 
$\log(1/\alpha_m)$. For the main applications that we will study in Section 
\ref{detection}, \ref{ellipsoides}, and \ref{simulations1}, the 
ratio $\log\left(1/\alpha_m\right) /\log\left(1/\alpha\right)$ is of 
order $\log(n)$, $\log\log n$, or $k\log(ep/k)$ where $k$ is a ``small'' integer. Thus, 
for each $\delta \in ]0,1-\alpha[$, the test based on $T_{\alpha}$ has a power 
greater than $1-\delta$ over a class of vectors which is close to $\bigcup_{m\in 
  \mathcal{M}}\mathcal{F}'_m(\delta,\alpha)$. It follows that for each $\theta\neq 0$ 
the power of this test under $\mathbb{P}_{\theta}$ is comparable to the best of 
the tests among the family $\left\{\phi_{m,\alpha},m\in \mathcal{M}\right\}$.

In the next two sections, we use this theorem to establish rates of testing 
against different types of alternatives. First, we give an upper bound for the 
rate of testing $\theta=0$ against a class of $\theta$ for which a lot of 
components are equal to $0$. In Section \ref{ellipsoides}, we study the 
rates of testing and simultaneous rates of testing $\theta=0$ against classes of ellipsoids. For the sake of simplicity, we will 
only consider the case $V=\{0\}$. Nevertheless, the procedure $T_{\alpha}$ 
defined in (\ref{definition_test}) applies in the same way when one considers more 
complex null hypothesis and the rates of testing are unchanged except that we 
have to replace $n$ by $n-d$ and $\var(Y)$ by $\var(Y|X_V)$. 
 
\section{Detecting non-zero coordinates}\label{detection}

Let us fix an integer $k$ between $1$ and $p$. In this section, we are 
interested in testing $\theta=0$ against the class of $\theta$ with a most $k$
non-zero components. This typically corresponds to the situation encountered when considering tests of neighborhood for large sparse graphs. As the graph is assumed to be sparse, only a small number of neighbors are missing under the alternative hypothesis. 

For each pair of integers $(k,p)$ with $k\leq p$, let $\mathcal{M}(k,p)$ be 
the class of all subsets of $\mathcal{I}=\{1,\ldots, p \}$ of cardinality $k$. The set 
$\Theta[k,p]$ stands for the subset of vectors $\theta \in 
\mathbb{R}^{p}$, such that at most $k$ coordinates of $\theta$ are 
non-zero. 
 
First, we define a test $T_{\alpha}$ of the form  (\ref{definition_test}) with Procedure $P_1$, 
and we derive an upper bound for the rate of testing of $T_{\alpha}$ against 
the alternative $\theta\in \Theta[k,p]$. Then, we show that this procedure is rate optimal when all the 
covariates are independent. Finally, we study the optimality of the test when 
$k=1$ for some examples of covariance matrix $\Sigma$.

\subsection{Rate of testing of $T_{\alpha}$}\label{rate_detection}

\begin{prte}\label{puissance_detection} 
We consider the set of models $\mathcal{M}=\mathcal{M}(k,p)$. We use the test 
$T_\alpha$ under Procedure $P_1$ and we take the weights $\alpha_m$ all 
equal to $\alpha/|\mathcal{M}|$. Let us suppose that $n$ satisfies: 
\begin{eqnarray}
n \geq L\left[\log\left(\frac{2}{\alpha\delta}\right)+k\log\left(\frac{ep}{k}\right) 
    \right] .\label{hypothese_detection}
\end{eqnarray}
Let us set the quantity 
\begin{eqnarray} 
  \rho'^2_{k,p,n} :=  L(\alpha,\delta)\frac{ k\log\left(\frac{ep}{k}\right) }{n}.\label{inegalite_puissance_detection} 
\end{eqnarray} 
For any $\theta$ in $\Theta[k,p]$, such that 
$\frac{\|\theta\|^2}{\emph{var}(Y)-\|\theta\|^2}\geq \rho'^2_{k,p,n}$, $\mathbb{P}_{\theta}\left(T_{\alpha}>0\right)\geq 1-\delta.$ 
\end{prte} 
We recall that the norm $\|.\|$ is defined in (\ref{definition_produit}) and  equals $\var(Y)-\var(Y|X)$.
This proposition easily follows from Theorem \ref{thrm_puissance} and its proof 
is given in Section \ref{proof1}. Note that the upper bound does not directly 
depend on the covariance matrix of the vector $X$. Besides, Hypothesis (\ref{hypothese_detection}) corresponds to the minimal assumption needed for consistency and type-oracle inequalities in the estimation setting as pointed out by Wainwright (\cite{wainwright07} Th. 2) and Giraud (\cite{giraud07} Sect. 3.1). Hence, we conjecture that Hypothesis (\ref{hypothese_detection}) is minimal so that Proposition \ref{puissance_detection} holds. We will further discuss the bound (\ref{inegalite_puissance_detection}) after deriving lower bounds for 
the minimax rate of testing. 
 
\subsection{Minimax lower bounds for independent covariates} \label{section_independent}
 
In the statistical framework considered here, the problem of giving minimax 
rates of testing under no prior knowledge of the covariance of $X$ and of 
$\var(Y)$ is open. This is why we shall only derive lower 
bounds when $\var(Y)$ and the covariance matrix of 
$X$ are known. In this section, we give non asymptotic lower bounds for 
the $(\alpha,\delta)$-minimax rate of testing over the set 
$\Theta[k,p]$ when the covariance matrix of $X$ is the identity matrix (except Proposition \ref{minoration_1test}). As 
these bounds coincide with the upper bound obtained in Section 
\ref{rate_detection}, this will show that our test $T_{\alpha}$ is rate optimal.~\\ 
 
We first give a lower bound for the $(\alpha,\delta)$-minimax rate of 
detection of all $p$ non-zero coordinates for any covariance matrix $\Sigma$. 
 
\begin{prte}\label{minoration_1test} 
Let us suppose that $\emph{var}(Y)$ is known.
Let us set $\rho^2_{p,n}$ such 
that: 
\begin{eqnarray} 
  \rho^2_{p,n} := L(\alpha,\delta)\frac{\sqrt{p}}{n}.\label{rho_1test} 
\end{eqnarray} 
Then for all $\rho < \rho_{p,n}$, 
$$\beta_\Sigma\left(\left\{\theta \in \Theta[p,p], \frac{\|\theta\|^2}{\emph{var}(Y)-\|\theta\|^2} = \rho^2 \right\}\right) \geq \delta,$$
where we recall that $\Sigma$ is the covariance matrix of $X$. 
\end{prte} 
 
If $n\geq (1+\gamma)p$ for some $\gamma>0$, Theorem 
   \ref{thrm_puissance} shows that the test $\phi_{\mathcal{I},\alpha}$ 
  defined in (\ref{definition_phi_test}) has power greater than 
  $\delta$ over the vectors $\theta$ that satisfy 
\begin{eqnarray*}
  \frac{\|\theta\|^2}{\var(Y)-\|\theta\|^2}\geq L(\gamma,\alpha,\delta)\frac{\sqrt{p}}{n}. 
\end{eqnarray*}  
Hence, $\sqrt{p}/n$ is the minimax rate of testing $\Theta[p,p]$ at least when the number of observations is larger than the number of covariates. This is coherent with the minimax rate obtained in the fixed design framework (e.g. \cite{baraudminimax}). When $p$ becomes larger we do not think that the lower bound given in Proposition \ref{minoration_1test} is still sharp. Note that this minimax rate of testing holds for any covariance matrix $\Sigma$ contrary to Theorem \ref{minoration_arete}.

We now turn to the lower bound for the $(\alpha,\delta)$-minimax rate of 
testing against $\theta\in \Theta[k,p]$. 
\begin{thrm}\label{minoration_arete} 
 Let us set $\rho^2_{k,p,n}$ such 
that 
\begin{eqnarray}\label{enonce_minoration} 
\rho^2_{k,p,n} := L(\alpha,\delta)\frac{ k}{n}  \log \left(1+ \frac{p}{k^2}+\sqrt{2\frac{p}{k^2}}\right). 
\end{eqnarray} 
We suppose that the covariance of $X$ is the identity matrix $I$. 
Then, for all $\rho < \rho_{k,p,n}$, 
$$\beta_I\left(\left\{\theta \in \Theta[k,p], \frac{\|\theta\|^2}{\emph{var}(Y)-\|\theta\|^2} = \rho^2 \right\}\right) > \delta.$$ 
where the quantity $\emph{var}(Y)$ is known.~\\~\\ 
If $\alpha + \delta \leq 53\%$, then one has 
$$\rho^2_{k,p,n} \geq \frac{k}{2n} \log\left(1+ \frac{p}{k^2} \vee \sqrt{\frac{p}{k^2}}\right).$$ 
\end{thrm} 
This result implies the following lower bound for the minimax rate of testing 
\begin{eqnarray*} 
\rho\left(\Theta[k,p],\alpha,\delta,\var(Y),I)\right)\geq \rho^2_{k,p,n}. 
\end{eqnarray*} 
The proof is given in Section \ref{proof2}. To the price of more technicalities, it is 
possible to prove that the lower bound still holds if the variables $(X_i)$ are 
independent with known variances possibly different. Theorem \ref{minoration_arete} recovers 
approximately the lower bounds for the minimax rates of testing in signal detection framework obtained by Baraud 
\cite{baraudminimax}. The main difference lies in the fact that 
we suppose $\var(Y)$ known which in the signal detection 
framework translates in the fact that we would know the quantity $\|f\|^2 
+\sigma^2$.

We are now in position to compare the results of Proposition 
\ref{puissance_detection} and Theorem \ref{minoration_arete}. Let distinguish between the 
values of $k$. 
\begin{itemize} 
\item When $k\leq p^{\gamma}$ for some $\gamma<1/2$, if $n$ is large enough to 
  satisfy the assumption of Proposition \ref{puissance_detection}, the quantities $\rho^2_{k,p,n}$ 
  and $\rho'^2_{k,p,n}$ are both of the order $\frac{k\log(p)}{n}$ times a constant 
  (which depends on $\gamma$, $\alpha$, and $\delta$). This shows that the 
  lower bound given in Theorem \ref{minoration_arete} is sharp. Additionally, 
  in this case, the procedure $T_{\alpha}$ defined in Proposition 
  \ref{puissance_detection} follows approximately the minimax rate of 
  testing. We recall that our procedure 
  $T_{\alpha}$ does not depend on the knowledge of $\var(Y)$ and 
  $\corr(X)$. In applications, a small $k$ typically 
  corresponds to testing a Gaussian graphical model with respect to a graph
  $\mathcal{G}$, when 
  the number of nodes is large and the graph is supposed to be sparse. When $n$ does not satisfy the assumption of Proposition \ref{puissance_detection}, we believe that our lower bound is not sharp anymore. 
 
\item When $\sqrt{p}\leq k \leq p$, the lower bound and the upper bound do not coincide 
   anymore. Nevertheless, if $n\geq (1+\gamma)p$ for some $\gamma>0$, Theorem 
   \ref{thrm_puissance} shows that the test $\phi_{\mathcal{I},\alpha}$ 
  defined in (\ref{definition_phi_test}) has power greater than 
  $\delta$ over the vectors $\theta$ that satisfy 
\begin{eqnarray}\label{majoration_globale} 
  \frac{\|\theta\|^2}{\var(Y)-\|\theta\|^2}\geq L(\gamma,\alpha,\delta)\frac{\sqrt{p}}{n}. 
\end{eqnarray}  
This upper bound and the lower bound do not depend on $k$. Here again, the 
lower bound obtained in Theorem \ref{minoration_arete} is sharp and the test 
$\phi_{\mathcal{I},\alpha}$ defined previously is rate optimal. The fact that the rate of testing 
stabilizes around $\sqrt{p}/n$ for $k>\sqrt{p}$ also appears in signal 
detection and there is a discussion of this phenomenon in 
\cite{baraudminimax}. 
\item When $k<\sqrt{p}$ and $k$ is close to $\sqrt{p}$, the lower bound and 
  the upper bound given by Proposition \ref{puissance_detection} differ from 
  at most a $\log(p)$ factor. For instance, if $k$ is of order 
  $\sqrt{p}/\log p$, the lower bound in Theorem \ref{minoration_arete} is 
  of order $\sqrt{p}\log\log p/\log p$ and the upper bound is of order 
  $\sqrt{p}$. We do not know if any of this bound is sharp and if the minimax 
  rates of testing coincide when $\var(Y)$ is fixed and when it is not fixed.   
\end{itemize} 
 
All in all, the minimax rates of testing exhibit the same range of rates in 
our framework as in signal detection 
\citep{baraudminimax} when the covariates are independent. Moreover, this implies 
that the minimax rate of testing is slower when the $(X_i)_{i\in \mathcal{I}}$ are independent than 
for any other form of dependence. Indeed, the upper bounds obtained in 
Proposition \ref{puissance_detection} and in (\ref{majoration_globale}) do not depend on the covariance of 
$X$. Then, a 
natural question arises: is the test statistic $T_{\alpha}$ rate optimal for 
other correlation of $X$? We will partially answer this question when 
testing against the alternative $\theta\in \Theta[1,p]$.

\subsection{Minimax rates for dependent covariates} 
In this section, we look for the minimax rate of testing $\theta=0$ against 
$\theta\in \Theta[1,p]$ when the covariates 
$X_i$ are no longer independent. 
We know that this rate is between the orders $\frac{1}{n}$, which is the 
minimax rate of testing when we know which coordinate is non-zero, and 
$\frac{\log(p)}{n}$, the minimax rate of testing for independent covariates. 
\begin{prte}\label{minoration_dependence} 
Let us suppose that there exists a positive number $c$ such that for any $i\neq 
j$, $$|\corr(X_i,X_j)|\leq c$$   and that $\alpha+\delta\leq 53\%$. 
We define $\rho_{1,p,n,c}^2$ as 
\begin{eqnarray} 
  \rho_{1,p,n,c}^2 := \frac{L}{n}\left(\log(p) \wedge 
  \frac{1}{c}\right).\label{minoration1_dependence} 
\end{eqnarray} 
Then for any $\rho < \rho_{1,p,n,c}$, 
$$\beta_{\Sigma}\left(\left\{\theta \in \Theta[1,p], 
\frac{\|\theta\|^2}{\text{\emph{var}}(Y)-\|\theta\|^2} = \rho^2 
\right\}\right) \geq \delta,$$ 
where $\Sigma$ refers to the covariance matrix of $X$. 
\end{prte}  
\textbf{Remark:} If the correlation between the covariates is smaller than 
$1/ \log(p)$, then the minimax rate of testing is of the same order as in the 
independent case. If the correlation between the covariates is larger, we 
show in the following Proposition that under some additional assumption, the 
rate is faster.

\begin{prte}\label{majoration_dependence} 
 
Let us suppose that the correlation between $X_i$ and $X_j$ is exactly $c>0$ for 
any $i\neq j$. Moreover, we assume that $n$ satisfies the following condition: 
\begin{eqnarray} 
n\geq L\left[1+\log \left(\frac{p}{\alpha\delta}\right)\right]\label{hypothese2majoration} 
\end{eqnarray} 
Let introduce the random variable $X_{p+1}:= \frac{1}{p}\sum_{i=1}^p\frac{X_i}{\sqrt{\var(X_i)}}$. If $\alpha<60\%$ and $\delta< 60\%$ the test $T_{\alpha}$ defined by 
$$T_{\alpha} =\left[\sup_{1\leq i\leq p} \phi_{\{i\},\alpha/(2p)}\right]\vee \phi_{\{p+1\},\alpha/2}$$ 
satisfies $$\mathbb{P}_0\left(T_{\alpha}>0\right)\leq \alpha \text{ and } 
\mathbb{P}_{\theta}\left(T_{\alpha}>0\right)\geq 1-\delta,$$ for any $\theta$ in $\Theta[1,p]$ such that 
$$\frac{\|\theta\|^2}{\text{\emph{var}}(Y)-\|\theta\|^2}\geq \frac{L(\alpha,\delta)}{n}\left(\log p\bigwedge\frac{1}{c}\right) .$$ 
\end{prte}

Consequently, when the correlation between $X_i$ and $X_j$ 
 is a positive constant $c$, the minimax rate of testing is 
 of order $\frac{\log(p)\wedge (1/c)}{n}$. When the correlation 
 coefficient $c$ is small, the minimax rate of testing coincides with the independent case, and when $c$ is larger 
 those rates differ. Therefore, the test 
 $T_{\alpha}$ defined in Proposition \ref{puissance_detection} is not rate optimal when the 
 correlation is known and is large. Indeed, when the correlation between the
 covariates is large, the tests statistics $\phi_{\{m\},\alpha_m}$ defining 
 $T_{\alpha}$ are highly correlated. The choice of the weights $\alpha_m$ in Procedure $P_1$ corresponds to 
 a Bonferroni procedure, which is precisely known to behave bad when 
 the tests are positively correlated. 
 
This example illustrates the limits of Procedure $P_1$. However, it 
is not very realistic to suppose that the covariates have a constant correlation, 
for instance when one considers a GGM. Indeed, we expect that the correlation between two covariates 
is large if they are neighbors in the graph and smaller if they are far (w.r.t. 
the graph distance). This is why we derive lower bounds of the rate of 
testing for other kind of correlation matrices often used to model stationary 
processes.

\begin{prte}\label{spatial} 
Let $X_1,\ldots, X_p$ form a stationary 
process on the one dimensional torus. More precisely, the correlation between 
$X_i$ and $X_j$ is a function of $|i-j|_p$ where $|.|_p$ refers to the 
toroidal distance defined by: 
$$|i-j|_p := (|i-j|) \wedge \left(p-|i-j|\right)\ .$$  
$\Sigma_1(w)$ and $\Sigma_2(t)$ respectively refer to the correlation matrix of $X$ such that 
\begin{eqnarray*} 
\text{\emph{corr}}(X_i,X_j)& =&  \exp(-w|i-j|_p) \text{ where }w>0\ ,\\ 
\text{\emph{corr}}(X_i,X_j)& = & (1+|i-j|_p)^{-t}\text{ where }t>0\ . 
\end{eqnarray*} 
 
Let us set $\rho^2_{1,p,n,\Sigma_1}(w)$ and $\rho^2_{1,p,n,\Sigma_2}(t)$  such that: 
\begin{eqnarray*} 
\rho^2_{1,p,n,\Sigma_1}(w)& :=& \frac{1}{n}\log\left(1+L(\alpha,\delta)p\frac{1-e^{-w}}{1+e^{-w}}\right)\\ 
\rho^2_{1,p,n,\Sigma_2}(t)& :=& \left\{\begin{array}{ccc}  
\frac{1}{n}\log\left(1+L(\alpha,\delta)\frac{p(t-1)}{t+1}\right)    & \text{if} &   t>1 \\ 
 \frac{1}{n}\log\left(1+L(\alpha,\delta)\frac{p}{1+2\log(p-1)}\right)    & \text{if} &    t=1\\ 
  \frac{1}{n}\log\left(1+L(\alpha,\delta)p^t2^{-t}(1-t)\right)   & \text{if} &    0<t<1. 
\end{array}\right. 
\end{eqnarray*} 
 
Then, for any $\rho^2< \rho^2_{1,p,n,\Sigma_1}(w)$, 
$$\beta_{\Sigma_1(w)}\left(\left\{\theta \in \Theta[1,p], \frac{\|\theta\|^2}{\emph{var}(Y)-\|\theta\|^2} = \rho^2 \right\}\right) \geq \delta,$$ 
and for any $\rho^2< \rho^2_{1,p,n,\Sigma_2}(t)$, 
$$\beta_{\Sigma_2(t)}\left(\left\{\theta \in \Theta[1,p], \frac{\|\theta\|^2}{\emph{var}(Y)-\|\theta\|^2} = \rho^2 \right\}\right) \geq \delta.$$ 
\end{prte} 
 
If the range $\omega$ is larger than $1/p^{\gamma}$  or if the range $t$ is larger than $\gamma$ for some $\gamma<1$, these lower bounds are of order $\frac{\log p}{n}$. As a consequence, for any of these correlation models the minimax rate of 
testing is of the same order as the minimax rate of testing for independent 
covariates. This means that our test $T_{\alpha}$ defined in Proposition \ref{puissance_detection} is rate-optimal for these correlations matrices. However, if $\omega$ is smaller than $1/p$ or if $t$ is smaller than $1/\log(p)$, we recover the parametric rates $1/n$, which is achieved by the test $\phi_{\{p+1\},\alpha}$. This comes from the fact that the correlation $\corr(X_1,X_i)$ does not converge to zero for such choices of $\omega$ or $t$. We omit the details since the arguments are similar to the proof of Proposition \ref{majoration_dependence}.~\\

To conclude, when $k\leq p^{\gamma}$ (for $\gamma\leq 1/2$), the test $T_{\alpha}$ defined in Proposition 
\ref{puissance_detection} is approximately $(\alpha,\delta)$-minimax against the 
alternative $\theta\in \Theta[k,p]$, when neither $\var(Y)$ nor the covariance 
matrix of $X$ is fixed. Indeed, the rate of testing of 
$T_{\alpha}$ coincide (up to a constant) with the supremum of the minimax rates of testing on $\Theta[k,p]$ over all possible covariance matrices $\Sigma$: 
$$\rho\left(\Theta[k,p],\alpha,\delta\right):= 
\sup_{\var(Y)>0,\Sigma>0}\rho\left(\Theta[k,p],\alpha,\delta,\var(Y),\Sigma \right),$$ 
where the supremum is taken over all positive $\var(Y)$ and every positive 
definite matrix $\Sigma$. When $k\geq \sqrt{p}$ and when $n\geq (1+\gamma)p$ 
(for $\gamma>0$), the test defined in (\ref{majoration_globale}) has the same behavior.

However, our procedure does not adapt to $\Sigma$: for some correlation 
matrices (as shown for instance in Proposition \ref{majoration_dependence}), 
$T_{\alpha}$ with Procedure $P_1$ is not rate optimal. Nevertheless, we 
believe and this will be illustrated in Section \ref{simulations1} that 
Procedure $P_2$ slightly improves the power of the test when the covariates are correlated.

\section{Rates of testing on ``ellipsoids" and adaptation}\label{ellipsoides} 
 
In this section, we define tests $T_{\alpha}$ of the form  (\ref{definition_test}) in order to 
test simultaneously $\theta=0$ against $\theta$ belongs to some classes of 
ellipsoids. We will study their rates and show that they are optimal at 
sometimes the price of a $\log p$ factor.
  
For any non increasing sequence $(a_i)_{1\leq i \leq {p+1}}$ such 
that $a_1=1$ and $a_{p+1}=0$ and any $R>0$, we define the ellipsoid $\mathcal{E}_{a}(R)$ by 
\begin{eqnarray}\label{definition_ellipsoides} 
 \mathcal{E}_{a}(R) := \left\{\theta\in \mathbb{R}^p,\sum_{i=1}^p \frac{\var(Y|X_{m_{i-1}}) - \var(Y|X_{m_i})}{a_i^2} \leq 
  R^2\var(Y|X)\right\}, 
\end{eqnarray} 
where $m_i$ refers to the set $\{1,\ldots,i\}$ and 
$m_0=\varnothing$. 
 
Let us explain why we call this set 
an ellipsoid. Assume for instance  
that the $(X_i)$ are independent identically distributed with variance one. In this case, the 
difference $\var(Y|X_{m_{i-1}}) - \var(Y|X_{m_i})$ equals $|\theta_i|^2$ 
and the definition of $\mathcal{E}_{a}(R)$ translates in 
$$ \mathcal{E}_{a}(R) = \left\{\theta\in \mathbb{R}^p,\sum_{i=1}^p 
  \frac{|\theta_i|^2}{a_i^2} \leq R^2 
  \var(Y|X)\right\}.$$ 
The main difference between this definition and the classical definition of an 
 ellipsoid in the fixed design regression framework (as for instance in 
 \cite{baraudminimax}) is the presence of the term  $\var(Y|X)$. We added this 
  quantity in order to be able to derive lower bounds of the minimax rate. If the $X_i$ are not i.i.d. with unit 
 variance, it is always possible to create a sequence $X'_i$ of 
 i.i.d. standard gaussian variables by orthogonalizing the $X_i$ using 
 Gram-Schmidt process. If we call $\theta'$ the vector in $\mathbb{R}^{p}$ such 
 that $X\theta=X'\theta'$, it is 
 straightforward to show that $\var(Y|X_{m_{i-1}}) - 
  \var(Y|X_{m_i})=|\theta'_i|^2$. We can then express $\mathcal{E}_a(R)$ using 
  the coordinates of $\theta'$ as previously: 
$$ \mathcal{E}_{a}(R) = \left\{\theta\in \mathbb{R}^p,\sum_{i=1}^p \frac{|\theta'_i|^2}{a_i^2} \leq 
 R^2 \var(Y|X)\right\}.$$ 
The main advantage of Definition \ref{definition_ellipsoides} is that it does not directly depend on the covariance of $X$. In the sequel we also consider the special case of ellipsoids with polynomial 
decay, 
\begin{eqnarray} 
\mathcal{E}'_{s}(R) := \left\{\theta\in \mathbb{R}^{p},\sum_{i=1}^p \frac{\var(Y|X_{m_{i-1}}) - \var(Y|X_{m_i})}{i^{-2s}\var(Y|X)} \leq 
  R^2\right\}, 
\end{eqnarray} 
where $s>0$ and $R>0$. First, we define two tests procedures of the form 
 (\ref{definition_test}) and evaluate 
their power respectively on the ellipsoids $\mathcal{E}_a(R)$ and on the 
ellipsoids $\mathcal{E}'_s(R)$. Then, we give some lower bounds for the 
$(\alpha,\delta)$-simultaneous minimax rates of testing. Extensions to more 
general $l_p$ balls with $0<p<2$ are possible to the price of more 
technicalities by adapting the results of Section 4 in Baraud \cite{baraudminimax}. 
 
These alternatives correspond to the situation where we are given an order of relevance on the covariates that are not in the null hypothesis. This order could either be provided by a previous knowledge of the model or by a model selection algorithm such as LARS (least angle regression) introduced by Efron \emph{et al.} \cite{LARS}. We apply this last method to build a collection of models for our testing procedure  (\ref{definition_test}) in \cite{villers_ggm}.

\subsection{Simultaneous Rates of testing of $T_{\alpha}$ over classes of ellipsoids}\label{rates_ellipsoides}

First, we define a procedure of the form  (\ref{definition_test}) in order to test 
$\theta=0$ against $\theta$ belongs to any of the ellipsoids 
$\mathcal{E}_a(R)$. For any $x>0$, $[x]$ denotes the integer part of $x$.

We choose the class of models $\mathcal{M}$ and the weights $\alpha_m$ as follows: 
\begin{itemize} 
\item If $n<2p$, we take the set $\mathcal{M}$ to be $\cup_{1\leq 
  k\leq[n/2]}m_k$ and all the weights $\alpha_m$ are equal to $\alpha/|\mathcal{M}|$. 
\item If $n\geq 2p$, we take the set $\mathcal{M}$ to be $\cup_{1\leq 
  k\leq p}m_k$. $\alpha_{m_p}$ equals $\alpha/2$ and for any $k$ between 1 and 
  $p-1$, $\alpha_{m_k}$ is chosen to be $\alpha/(2(p-1))$. 
\end{itemize} 
 
As previously, we bound the power of the tests $T_{\alpha}$ from a non-asymptotic point of view.
 
\begin{prte}\label{puissance_ellipsoide} 
Let us assume that 
\begin{eqnarray}\label{condition_n_ellipsoide} 
n \geq L \left[1+ \log\left(\frac{1}{\alpha\delta}\right)\right].
\end{eqnarray} 
For any ellipsoid $\mathcal{E}_a(R)$, the test $T_{\alpha}$ defined by  (\ref{definition_test}) with Procedure 
$P_1$ and with the class of models given just above satisfies  
$$\mathbb{P}_{0}\left(T_{\alpha}\leq 0\right)\geq 1- \alpha,$$ 
and $\mathbb{P}_{\theta}\left(T_{\alpha}>0\right)\geq 1-\delta$ for all 
$\theta \in \mathcal{E}_a(R)$ such that 
\begin{equation} 
\frac{\|\theta\|^2}{\emph{var}(Y)-\|\theta\|^2} \geq L(\alpha,\delta)\log n  
\inf_{1\leq i \leq [n/2]}\left[a_{i+1}^2R^2 + 
  \frac{\sqrt{i}}{n}\right]  
\end{equation} 
if $n<2p$, or  
\begin{equation}\label{majoration_ngrand} 
 \frac{\|\theta\|^2}{\emph{var}(Y)-\|\theta\|^2} \geq L(\alpha,\delta)\left\{\left[
\log p\inf_{1\leq i \leq p-1}\left(a_{i+1}^2R^2 + 
  \frac{\sqrt{i}}{n}\right)\right] \bigwedge \frac{\sqrt{p}}{n}\right\} 
\end{equation} 
if $n\geq 2p$.  
\end{prte} 
 
All in all, for large values of $n$, the rate of testing is of order $\sup_{1\leq i \leq 
  p}\left[a_i^2R^2 \wedge 
  \frac{\sqrt{i\log(p)}}{n}\right]$.  We show in the next subsection that 
  the minimax rate of testing for an ellipsoid is of order: 
$$\sup_{1\leq i \leq {p}}\left[a_{i}^2R^2 \wedge \frac{\sqrt{i}}{n}\right].$$ 
Besides, we prove in Proposition \ref{minoration_minimax_collection_ellipsoide} that a loss in $\sqrt{\log\log p}$ is unavoidable if one considers the simultaneous minimax rates of testing over a family of nested ellipsoids. Nevertheless, we 
do not know if the term $\sqrt{\log(p)}$ is optimal for testing simultaneously 
against all the ellipsoids $\mathcal{E}_a(R)$ for all sequences 
$(a_i)$ and all $R>0$. When $n$ is smaller than $2p$, we obtain 
comparable results except that we are unable to consider alternatives in large dimensions in the infimum (\ref{majoration_ngrand}).~\\

We now turn to define a procedure of the form  (\ref{definition_test}) in order to test simultaneously that 
$\theta=0$ against $\theta$ belongs to any of the $\mathcal{E}'_s(R)$. For 
this, we introduce the following collection of models $\mathcal{M}$ and 
weights $\alpha_m$: 
\begin{itemize} 
\item If $n<2p$, we take the set $\mathcal{M}$ to be $\cup m_k$ where $k$ 
  belongs to $ \left\{2^j,j\geq 0\right\}\cap \left\{1,\ldots 
,[n/2]\right\}$ and all the weights $\alpha_m$ are chosen to be $\alpha/|\mathcal{M}|$. 
\item If $n\geq 2p$, we take the set $\mathcal{M}$ to be $\cup m_k$ where $k$ 
  belongs to $ \left(\left\{2^j,j\geq 0\right\}\cap \left\{1,\ldots 
,p\right\}\right)\cup \{p\}$, $\alpha_{m_p}$ equals $\alpha/2$ and for any $k$ 
  in the model between 1 and 
  $p-1$, $\alpha_{m_k}$ is chosen to be $\alpha/(2(|\mathcal{M}|-1))$. 
\end{itemize} 
 
\begin{prte}\label{puissance_boules}  
Let us assume that 
\begin{eqnarray} \label{condition_n_boules}
n \geq L \left[1+\log\left(\frac{1}{\alpha\delta}\right)\right]
\end{eqnarray} 
and that $R^2\geq  \sqrt{\log\log n}/n$. 
For any $s>0$, the test procedure $T_{\alpha}$ defined by  (\ref{definition_test}) with Procedure 
$P_1$ and with a class of models given just above satisfies: 
$$\mathbb{P}_{0}\left(T_{\alpha}>0\right)\geq 1- \alpha,$$ 
and $\mathbb{P}_{\theta}\left(T_{\alpha}>0\right)\geq 1-\delta$ for any 
$\theta\in \mathcal{E}'_s(R)$ such that 
\begin{equation}\label{majoration_npetit_boule} 
\frac{\|\theta\|^2}{\emph{var}(Y)-\|\theta\|^2} \geq 
L\left(\alpha,\delta\right)\left[R^{2/(1+4s)} \left(\frac{\sqrt{\log\log n 
     }}{n}\right)^{4s/(1+4s)}+ R^2 \left(n/2\right)^{-2s}+\frac{\log\log n}{n}\right] 
\end{equation} 
if $n<2p$ or 
\begin{equation}\label{majoration_ngrand_boule} 
\frac{\|\theta\|^2}{\emph{var}(Y)-\|\theta\|^2} \geq 
L\left(\alpha,\delta\right)\left(\left[R^{2/(1+4s)} \left(\frac{\sqrt{\log\log 
      p}}{n}\right)^{4s/(1+4s)}+\frac{\log\log p}{n}\right]\bigwedge 
\frac{\sqrt{p}}{n}\right) 
\end{equation} 
if $n\geq 2p$.  
\end{prte}

Again, we retrieve similar results to those of Corollary 2 in 
\cite{baraud2003} in the fixed design regression framework. For $s>1/4$ and $n< 2p$, the rate of testing is of order $\left(\frac{\sqrt{\log\log 
 n} }{n}\right)^{4s/(1+4s)}$. We show in the next subsection that the
logarithmic factor is due to the adaptive property of the 
 test. If $s\leq 1/4$, the rate is of order $n^{-2s}$. When $n\geq 2p$,  
 the rate is of order $\left(\frac{\sqrt{\log\log 
 p}}{n}\right)^{4s/(1+4s)} \wedge \left(\frac{\sqrt{p}}{n}\right)$, and we mention at 
 the end of the next subsection that it is optimal. 
 
Here again, it is possible to define these tests with Procedure $P_2$ in 
order to improve the power of the test (see Section \ref{simulations1} 
for numerical results).

\subsection{Minimax lower bounds} 
 
We first establish the $(\alpha,\delta)$-minimax  rate of testing over an 
ellipsoid when the variance of $Y$ and the covariance matrix of $X$ are known.

\begin{prte}\label{vitesse_minimax_ellipsoide} 
Let us set the sequence $(a_i)_{1\leq i\leq {p+1}}$ and the positive number $R$.  
We introduce 
\begin{eqnarray} 
\rho^2_{a,n}(R) := \sup_{1\leq i\leq p}\left[\rho^2_{i,n}\wedge a_i^2R^2\right], 
\end{eqnarray} 
where $\rho^2_{i,n}$ is defined by (\ref{rho_1test}), then for any non singular covariance matrix $\Sigma$ we have 
\begin{eqnarray*} 
  \beta_{\Sigma}\left(\left\{\theta \in 
  \mathcal{E}_a(R),\frac{\|\theta\|^2}{\emph{var}(Y)- 
  \|\theta\|^2}\geq \rho_{a,n}^2(R)\right\}\right)\geq \delta, 
\end{eqnarray*} 
where the quantity $\var(Y)$ is fixed. If $\alpha +\delta\leq 47 \%$ then 
\begin{eqnarray*} 
  \rho^2_{a,n}(R) \geq \sup_{1\leq i\leq p}\left[\frac{\sqrt{i}}{n}\wedge a_i^2R^2\right]. 
\end{eqnarray*} 
\end{prte}

This lower bound is once more analogous to the one in the fixed design 
regression framework.  Contrary to the lower bounds obtained in the previous section, it does 
not depend on the covariance of the covariates. We now look for an  upper bound of the 
minimax rate of testing over a given ellipsoid. First, we need to define the quantity $D^*$ as: 
\begin{eqnarray*} 
D^* := \inf\left\{1\leq i \leq p, a_i^2R^2\leq \frac{\sqrt{i}}{n} \right\} 
\end{eqnarray*} 
with the convention that $\inf \varnothing=p$.

\begin{prte}\label{majoration_minimax_ellipsoide} 
Let us assume that $n\geq L\log\left[1+\log\left(\frac{1}{\alpha\delta}\right)\right]$. 
If $R^2 > \frac{1}{n}$ and $D^*\leq n/2$, the test $\phi_{m_{D^*},\alpha}$ defined by (\ref{definition_phi_test}) 
satisfies 
\begin{eqnarray*} 
\mathbb{P}_{0}\left[\phi_{m_{D^*},\alpha}= 1\right]\leq \alpha \text{ and } 
\mathbb{P}_{\theta}\left[\phi_{m_{D^*},\alpha} = 0\right]\leq \delta 
\end{eqnarray*} 
for all $\theta\in \mathcal{E}_a(R)$ such that 
\begin{eqnarray*} 
  \frac{\|\theta\|^2}{\emph{var}(Y)-\|\theta\|^2}\geq L(\alpha,\delta) \sup_{1\leq i\leq 
  p}\left[\frac{\sqrt{i}}{n}\wedge a_i^2R\right]. 
\end{eqnarray*} 
\end{prte} 
 
If $n\geq 2D^*$, the rates 
of testing on an ellipsoid are analogous to the rates on an ellipsoid in fixed design 
regression framework (see for instance \cite{baraudminimax}). If $D^*$ is 
large and $n$ is small, the bounds in Proposition 
\ref{vitesse_minimax_ellipsoide} and \ref{majoration_minimax_ellipsoide} do 
not coincide. In this case, we do not know if this comes from the fact that 
the test in Proposition \ref{majoration_minimax_ellipsoide} does not 
depend on the knowledge of $\var(Y)$ or if one of the bounds in Proposition 
\ref{vitesse_minimax_ellipsoide} and \ref{majoration_minimax_ellipsoide} is not sharp.

We are now interested in computing lower bounds for rates of testing 
simultaneously over a family of ellipsoids, in order to compare them with 
rates obtained in Section \ref{rates_ellipsoides}. First, we need a lower bound for the minimax 
simultaneous rate of testing over nested linear spaces. We recall that for any 
$D\in \{1,\ldots, p\}$, $S_{m_D}$ stands for the linear spaces of vectors $\theta$ 
such that only their $D$ first coordinates are possibly non-zero. 
 
\begin{prte}\label{minimax_lower_nested_linear} 
For $D\geq 2$, let us set  
\begin{eqnarray} 
\bar{\rho}^2_{D,n} :=  L(\alpha,\delta)\frac{\sqrt{\log \log (D+1) }\sqrt{D}}{n}.\label{definition_rho_bar} 
\end{eqnarray} 
Then, the following lower bound holds 
\begin{eqnarray*} 
  \beta_{I}\left(\bigcup_{1\leq D\leq p} \left\{\theta \in S_{m_D}, \frac{\|\theta\|^2}{\emph{var}(Y)-\|\theta\|^2}=r_D^2\right\}\right)\geq \delta, 
\end{eqnarray*} 
if for all $D$ between $1$ and $p$, $r_D\leq \bar{\rho}_{D,n}$ 
\end{prte} 
 
Using this Proposition, it is possible to get a lower bound for the simultaneous rate of 
testing over a family of nested ellipsoids. 
 
\begin{prte}\label{minoration_minimax_collection_ellipsoide} 
We fix a sequence $(a_i)_{1\leq i\leq p+1}$. For each $R>0$, let us set 
\begin{eqnarray}\label{definition_rho_ellipsoide} 
\bar{\rho}^2_{a,R,n} := \sup_{1\leq D\leq p}\left[\bar{\rho}^2_{D,n} \wedge (R^2a_D^2)\right]. 
\end{eqnarray} 
where $\bar{\rho}_{D,n}$ is given by (\ref{definition_rho_bar}). Then, for any non singular 
covariance matrix $\Sigma$ of the vector $X$,  
\begin{eqnarray*} 
\beta_{\Sigma}\left(\bigcup_{R>0}\left\{\theta\in \mathcal{E}_a(R),\frac{\|\theta\|^2}{\emph{var}(Y)-\|\theta\|^2}\leq \bar{\rho}^2_{a,R,n}\right\}\right)\geq \delta. 
\end{eqnarray*} 
\end{prte} 
 
This Proposition shows that the problem of adaptation is impossible in this 
setting: it is impossible to define a test which is simultaneously minimax 
over a class of nested ellipsoids (for $R>0$). This is also the case in  
fixed design as proved by \cite{spokoiny96} for the case 
of Besov bodies. The loss of a term of the order $\sqrt{\log\log p}/n$ is unavoidable.

As a special case of Proposition \ref{minoration_minimax_collection_ellipsoide}, it is possible to compute a 
lower bound for the simultaneous minimax rate over $\mathcal{E}'_s(R)$ where 
$R$ describes the positive numbers. After some calculation, we find 
a lower bound of order: 
$$\left(\frac{\sqrt{\log\log p} }{n}\right)^{\frac{4s}{1+4s}}\bigwedge \frac{\sqrt{p\log\log p}}{n}.$$ 
This shows that the power of the test 
$T_{\alpha}$ obtained in (\ref{majoration_ngrand_boule}) for 
$n\geq 2p$ is optimal when 
$R^2\geq \sqrt{\log\log n}/n$. However, when $n<2p$ and $s\leq 1/4$, we do not 
know if the rate $n^{-2s}$ is optimal or not.  
 
To conclude, when $n\geq 2p$ the test $T_{\alpha}$ defined in 
Proposition \ref{puissance_boules} achieves the simultaneous minimax rate over the classes of 
ellipsoids $\mathcal{E}'_s(R)$. On the other hand, the test $T_{\alpha}$ 
defined in Proposition \ref{puissance_ellipsoide} is not rate optimal 
simultaneously over all the ellipsoids $\mathcal{E}_a(R)$ and 
suffers a loss of a $\sqrt{\log p}$ factor even when $n\geq 2p$.

\section{Simulations studies}\label{simulations1} 
 
The purpose of this simulation study is threefold. First, we illustrate the theoretical 
results established in previous sections. Second, we show that our 
procedure is easy to implement for different  
choices of collections $\mathcal{M}$ and is computationally feasible even when $p$ is large. Our third purpose is to compare the 
efficiency of Procedures $P_1$  and $P_2$. Indeed, for a given collection 
$\mathcal{M}$, we know from Section  
  \ref{compa_puissance} that the test (\ref{definition_test}) based 
  on Procedure $P_2$ is more powerful than the corresponding test based on 
  $P_1$. However, the computation  
of the quantity $q_{{\bf X}, \alpha}$ is possibly time consuming and we
therefore want to know if the benefit in power is worth the computational burden. 

To our knowledge, when the number of covariates $p$ is larger than the 
number of observations $n$ there is no test with which we can compare 
our procedure.

\subsection{Simulation experiments}

We consider the regression model (\ref{model}) with $\mathcal{I}=\{1, \ldots, 
p\}$ and test the null hypothesis  "$\theta=0$", which is equivalent 
to "$Y$ is independent of $X$", at  
level $\alpha= 5\%$. Let  $(X_i)_{1 \leqslant i \leqslant p}$ be a 
collection of $p$ Gaussian variables with unit variance.  The random 
variable is defined as  
follows: $Y=\sum_{i=1}^p\theta_i X_i + \varepsilon $ where $\varepsilon$ 
is a zero mean gaussian variable with variance $1-\|\theta\|^2$ 
independent of $X$.  
 
We consider two simulation experiments described below. 
\begin{enumerate} 
\item First simulation experiment: The  correlation 
between $X_i$ and $X_j$ is a constant $c$ for any $i\neq j$.  
Besides, in this experiment the parameter  $\theta$ is chosen such that 
only one of its components is possibly non-zero. This corresponds to 
the situation considered in Section \ref{detection}. 
First, the number of covariates $p$ is fixed equal to $30$ and the 
number of observations $n$ is taken equal to  $10$ and $15$.  
 We choose for $c$ three different values $0$, 
$0.1$, and $0.8$, allowing thus to compare the procedure for independent, 
weakly and highly correlated covariates. We estimate the size of the 
test by taking $\theta_1=0$ and the power by taking for $\theta_1$ the 
values $0.8$ and $0.9$. Theses  choices of $\theta$ lead to  a 
small and a large signal/noise ratio $r_{s/n}$ defined in 
(\ref{ratio}) and equal in this experiment to  $\theta_1^2/  (1-\theta_1^2)$. 
Second, we examine the behavior of the tests  
when $p$ increases and when the covariates are highly correlated: $p$ 
equals $100$ and $500$, $n$ equals  $10$ and $15$, $\theta_1$ is set to 
$0$ and $0.8$, and $c$ is chosen to be $0.8$. 
\item Second simulation experiment: The covariates 
  $(X_i)_{1 \leqslant i \leqslant p}$ are independent. The 
  number of covariates $p$ equals $500$ and the number of  
observations $n$ equals $50$ and $100$. We set  for any $i \in \{1, 
\ldots, p\}$,  
$\theta_i=R i^{-s}$. We estimate the size of the test by taking $R=0$ 
and the power by taking for $(R,s)$ the value $(0.2,0.5)$, which 
corresponds to a slow decrease of the $(\theta_i)_{1 \leqslant i 
  \leqslant p}$. It was pointed out in the beginning of Section 
\ref{ellipsoides} that $|\theta_i|^2$ equals 
$\var(Y|X_{m_{i-1}}) - \var(Y|X_{m_i})$. Thus, $|\theta_i|^2$ represents the benefit in term of conditional variance brought by the 
variable $X_i$.  
\end{enumerate}

We use our testing procedure defined in (\ref{definition_test}) with different 
collections $\mathcal{M}$ and different choices for the weights 
$\{\alpha_m,m\in \mathcal{M}\}$. \vspace{0.5cm}
 
\textit{The collections $\mathcal{M}$:} we define three classes. 
Let us set $J_{n,p} = p\wedge  [\frac{n}{2}]$, where 
$[x]$ denotes the integer part of $x$ and let us define:  
\begin{eqnarray*} 
\mathcal{M}^1 & :=& \{\{i\}, 1\leqslant i  \leqslant p\}\\ 
\mathcal{M}^2 & := & \{m_k=\{1, 2, \ldots, k\},   1\leqslant k  \leqslant 
J_{n,p}\}\}\\  
\mathcal{M}^3 & := & \{m_k=\{1, 2, \ldots, k\},   k \in \{2^j, j \ge 0\} \cap \{1, 
\ldots, J_{n,p}\}  \} 
\end{eqnarray*} 
We evaluate the performance of our testing procedure with 
$\mathcal{M}=\mathcal{M}^1$ in  
the first simulation experiment, and $\mathcal{M}=\mathcal{M}^2$ and 
$\mathcal{M}^3$ in the  
second simulation experiment. The cardinality of these three collections is smaller than $p$, and the computational complexity of the testing procedures is at most linear in $p$.\vspace{0.5cm} 
 
\textit{The collections $\{ \alpha_m, m \in \mathcal{M} \}$:} We 
consider Procedures $P_1$ and $P_2$ defined in Section 
\ref{section_definition}. When we are  
using the procedure $P_1$, the $\alpha_m$'s equal $\alpha/|\mathcal{M}|$ where 
$|\mathcal{M}|$ denotes the cardinality of the collection $\mathcal{M}$ 
. The quantity  
$q_{{\bf X}, \alpha}$ that occurs in the procedure $P_2$ is computed by 
simulation. We use $1000$ simulations for the estimation of $q_{{\bf X}, 
  \alpha}$. In the sequel we note $T_{\mathcal{M}^i,P_j}$ the test 
(\ref{definition_test}) with collection $\mathcal{M}^i$ and Procedure $ P_j$.

In the first experiment, when $p$ is large we also consider two other tests: 
\begin{enumerate} 
\item The test $\phi_{\{1\},\alpha}$ (\ref{definition_phi_test}) of the 
  hypothesis $\theta_1=0$ against the alternative $\theta_1 \neq 
0$. This test corresponds to the single test when we know which 
coordinate is non-zero.   
\item The test $\phi_{\{p+1\},\alpha}$ where $X_{p+1}:= \frac{1}{p}\sum_{i=1}^p X_i$. Adapting the proof of Proposition 
\ref{majoration_dependence}, we know that this test is approximately minimax 
on $\Theta[1,p]$ if the correlation between the covariates is constant and 
large. 
\end{enumerate} 
Contrary to our procedures, these two tests are based on the  
knowledge of $\var(X)$ (and eventually $\theta$). We only use them as a benchmark to  
evaluate the performance of our procedure. We aim at showing that our test with 
Procedure $P_2$ is as powerful than $\phi_{\{p+1\},\alpha}$ and is close 
to the test $\phi_{\{1\},\alpha}$.

We  estimate the size and the power of the testing procedures with 
$1000$ simulations. For each simulation, we simulate the gaussian vector 
$(X_1, \ldots, X_p)$ 
and then simulate the variable $Y$ as described in the two simulation 
experiments.

\subsection{Results of the simulation}  
 
\begin{Table}[h] 
\caption{ First simulation study, independent case: 
  $p=30$,  $c=0$. Percentages of rejection and value of the signal/noise ratio 
  $r_{s/n}$.\label{tab1}}

\begin{center} 
Null hypothesis is true, $\theta_1=0$ 
\end{center} 
\begin{center} 
\begin{tabular}{l|rr} 
$n$ & $T_{\mathcal{M}^1,P_1}$  & $T_{\mathcal{M}^1,P_2}$\\ 
\hline 10 &  0.043& 0.045  \\ 
15 & 0.044 & 0.049   \\ 
\end{tabular} 
\end{center}

\begin{center} 
Null hypothesis is false 
\end{center} 
\begin{minipage}{8cm} 
\begin{center} 
\begin{tabular}{l|rr} 
\multicolumn{3}{c}{ $\theta_1=0.8$, $r_{s/n}=1.78$}\\ 
$n$ & $T_{\mathcal{M}^1,P_1}$  & $T_{\mathcal{M}^1,P_2}$\\ 
\hline 10 & 0.48  & 0.48  \\ 
15 & 0.81 &  0.81 \\ 
\end{tabular} 
\end{center} 
\end{minipage} 
\ ~ \ 
\begin{minipage}{8cm} 
\begin{center} 
\begin{tabular}{l|rr} 
\multicolumn{3}{c}{$\theta_1=0.9$, $r_{s/n}=4.26$}\\ 
$n$ & $T_{\mathcal{M}^1,P_1}$  & $T_{\mathcal{M}^1,P_2}$\\ 
\hline 10 & 0.86 & 0.86  \\ 
15 & 0.99 &  0.99\\ 
\end{tabular} 
\end{center} 
\end{minipage} 
\end{Table}

The results of the first simulation experiment for $c=0$ are given in 
Table  \ref{tab1}. As expected, the power of the tests increases with 
the number  
of observations $n$ and with the signal/noise ratio $r_{s/n}$. 
If the signal/noise ratio is large enough, we obtain powerful tests even 
if the number of covariates $p$ is larger than the number of observations.

\begin{Table}[hptb] 
\caption{ First simulation study, independent and dependent case. $p=30$ 
  Frequencies of rejection. \label{tab2}}   
  
\begin{center} 
Null hypothesis is true, $\theta_1=0$ 
\end{center} 
 
\begin{minipage}{5cm} 
\begin{center} 
\begin{tabular}{l|rr} 
\multicolumn{3}{c}{ $c=0$ }\\ 
$n$ & $T_{\mathcal{M}^1,P_1}$  & $T_{\mathcal{M}^1,P_2}$\\ 
\hline 10 &  0.043& 0.045  \\ 
15 & 0.044 & 0.049  \\ 
\end{tabular} 
\end{center} 
\end{minipage} 
\ ~ \ 
\begin{minipage}{5cm} 
\begin{center} 
\begin{tabular}{l|rr} 
\multicolumn{3}{c}{ $c=0.1$ }\\ 
$n$ & $T_{\mathcal{M}^1,P_1}$  & $T_{\mathcal{M}^1,P_2}$\\ 
\hline 10  & 0.042  &0.04\\ 
15 &0.058  & 0.06   \\ 
\end{tabular} 
\end{center} 
\end{minipage} 
\ ~ \ 
\begin{minipage}{5cm} 
\begin{center} 
\begin{tabular}{l|rr} 
\multicolumn{3}{c}{ $c=0.8$}\\ 
$n$ & $T_{\mathcal{M}^1,P_1}$  & $T_{\mathcal{M}^1,P_2}$\\ 
\hline 10 &0.018  &0.045  \\ 
15 &0.019 & 0.052  \\ 
\end{tabular} 
\end{center} 
\end{minipage}

\begin{center} 
Null hypothesis is false, $\theta_1=0.8$ 
\end{center} 
 
\begin{minipage}{5cm} 
\begin{center} 
\begin{tabular}{l|rr} 
\multicolumn{3}{c}{ $c=0$ }\\ 
$n$ & $T_{\mathcal{M}^1,P_1}$  & $T_{\mathcal{M}^1,P_2}$\\ 
\hline 10  & 0.48 &0.48  \\ 
15 & 0.81 &  0.81  \\ 
\end{tabular} 
\end{center} 
\end{minipage} 
\ ~ \ 
\begin{minipage}{5cm} 
\begin{center} 
\begin{tabular}{l|rr} 
\multicolumn{3}{c}{ $c=0.1$ }\\ 
$n$ & $T_{\mathcal{M}^1,P_1}$  & $T_{\mathcal{M}^1,P_2}$\\ 
\hline 10  & 0.49 &0.49 \\ 
15 & 0.81&  0.82\\ 
\end{tabular} 
\end{center} 
\end{minipage} 
\ ~ \ 
\begin{minipage}{5cm} 
\begin{center} 
\begin{tabular}{l|rr} 
\multicolumn{3}{c}{ $c=0.8$}\\ 
$n$ & $T_{\mathcal{M}^1,P_1}$  & $T_{\mathcal{M}^1,P_2}$\\ 
\hline 10 & 0.64 &0.77  \\ 
15 & 0.89&  0.94\\ 
\end{tabular} 
\end{center} 
\end{minipage} 
\end{Table}

In Table \ref{tab2} we present results of the first simulation experiment for 
$\theta_1=0.8$ when $c$ varies. \\ 
Let us first compare the results for 
independent, weakly and highly correlated covariates when using 
Procedure $P_1$.  
The size and the power of 
the test for weakly correlated covariates are similar to the size and 
the power obtained in the independent case.  
Hence, we recover the remark following Proposition 
\ref{minoration_dependence}: when the  
correlation coefficient between the covariates is small, the minimax rate 
is of the same order as in  the independent case. 
The test for highly correlated covariates is 
more powerful than the test for independent covariates, recovering thus the 
remark following Theorem \ref{minoration_arete}: the worst case from a minimax 
rate perspective is the case where the covariates are independent.
Let us now compare Procedures $P_1$ and $P_2$.  In the case of 
independent or weakly correlated covariates, they give similar 
results. For highly correlated covariates,   
the power of 
$T_{\mathcal{M}^1,P_2}$ is much larger than the one of $T_{\mathcal{M}^1,P_1}$.

\begin{Table}[hptb] 
\caption{First simulation study, dependent case: $ c=0.8$.  Frequencies 
  of rejection. \label{tab3}}  
 
\begin{center} 
Null hypothesis is true, $\theta_1=0$ 
\end{center} 
 
\begin{minipage}{8cm} 
\begin{center} 
\begin{tabular}{l|rrrr} 
\multicolumn{5}{c}{ $p=100$ }\\ 
$n$ & $T_{\mathcal{M}^1,P_1}$  & $T_{\mathcal{M}^1,P_2}$& 
$\phi_{\{1\},\alpha}$ &$\phi_{\{p+1\},\alpha}$ \\  
\hline 10  & 0.01  &0.056 &0.051 &0.045\\ 
15& 0.016& 0.053& 0.047&0.053\\ 
\end{tabular} 
\end{center} 
\end{minipage} 
\ ~ \ 
\begin{minipage}{8cm} 
\begin{center} 
\begin{tabular}{l|rrrr} 
\multicolumn{5}{c}{ $p=500$}\\ 
$n$ & $T_{\mathcal{M}^1,P_1}$  & $T_{\mathcal{M}^1,P_2}$ & 
$\phi_{\{1\},\alpha}$ & $\phi_{\{p+1\},\alpha}$ \\  
\hline 10 & 0.009  &0.044 & 0.040& 0.040\\ 
15 &  0.011& 0.040 &0.042&0.034\\ 
\end{tabular} 
\end{center} 
\end{minipage} 
 
\begin{center} 
Null hypothesis is false, $\theta_1=0.8$ 
\end{center} 
 
\begin{minipage}{8cm} 
\begin{center} 
\begin{tabular}{l|rrrr} 
\multicolumn{5}{c}{ $p=100$ }\\ 
$n$ & $T_{\mathcal{M}^1,P_1}$  & $T_{\mathcal{M}^1,P_2}$& 
$\phi_{\{1\},\alpha}$ &$\phi_{\{p+1\},\alpha}$ \\  
\hline 10  & 0.60  &0.77 &0.91 &0.79\\ 
15& 0.85& 0.92& 0.99&0.92\\ 
\end{tabular} 
\end{center} 
\end{minipage} 
\ ~ \ 
\begin{minipage}{8cm} 
\begin{center} 
\begin{tabular}{l|rrrr} 
\multicolumn{5}{c}{ $p=500$}\\ 
$n$ & $T_{\mathcal{M}^1,P_1}$  & $T_{\mathcal{M}^1,P_2}$ & 
$\phi_{\{1\},\alpha}$ & $\phi_{\{p+1\},\alpha}$ \\  
\hline 10 & 0.52 &0.76 &0.91 &0.77 \\ 
15 & 0.77 &  0.94 &0.99&0.94\\ 
\end{tabular} 
\end{center} 
\end{minipage} 
\end{Table}

In Table \ref{tab3} we present results of the multiple testing 
procedure and of the two tests $\phi_{\{1\},\alpha}$ and 
$\phi_{\{p+1\},\alpha}$ when  
$c=0.8$ and the number of covariates $p$ is large. For $p=500$ and $n=15$, one test takes less than one second with Procedure $P_1$ and less than $30$ seconds with Procedure $P_2$. As expected, 
Procedure $P_1$ is too conservative when $p$ increases. 
For $p=100$, the power of the test based on Procedure $P_1$ is smaller than the power of the test $\phi_{\{p+1\},\alpha}$ and this difference increases  
when $p$ is larger. The  test based on  Procedure 
$P_2$ is as powerful as $\phi_{\{p+1\},\alpha}$, and its power 
is close to the one of $\phi_{\{1\},\alpha}$. We recall that this last test 
is based on the knowledge of the non-zero component of $\theta$ contrary to 
ours. Besides, the test $\phi_{\{p+1\},\alpha}$ was shown in Proposition \ref{majoration_dependence} to be optimal for this particular correlation setting. Hence, Procedure $P_2$ seems to achieve the optimal rate in this situation. Thus, we advise to use in practice Procedure $P_2$ if the number of 
covariates $p$ is large, because Procedure $P_1$ becomes too conservative, 
especially if the covariates are correlated.

\begin{Table}[h] 
\caption{Second simulation study. Frequencies of rejection. \label{tab4}}  
 
\begin{center} 
Null hypothesis is true, $R=0$ 
\end{center} 
 
\begin{center} 
\begin{tabular}{l|rrrr} 
$n$ & $T_{\mathcal{M}^2,P_1}$  & $T_{\mathcal{M}^2,P_2}$ & 
$T_{\mathcal{M}^3,P_1}$ & $T_{\mathcal{M}^3,P_2}$\\  
\hline 50  &  0.013 &0.052 & 0.036 & 0.059 \\ 
 100  &0.009   & 0.059   & 0.042 & 0.059\\ 
  
\end{tabular} 
\end{center}

\begin{center} 
Null hypothesis is false, $R=0.2, s=0.5$ 
\end{center} 
 
\begin{center} 
\begin{tabular}{l|rrrr} 
$n$ & $T_{\mathcal{M}^2,P_1}$  & $T_{\mathcal{M}^2,P_2}$ & 
$T_{\mathcal{M}^3,P_1}$ & $T_{\mathcal{M}^3,P_2}$\\  
\hline 50  &0.17   & 0.33 & 0.31 & 0.38\\ 
 100  &0.42   & 0.66 & 0.62 & 0.69\\ 
 
\end{tabular} 
\end{center} 
\end{Table} 
 
The results of the second simulation experiment are given in 
Table \ref{tab4}. 
As expected, Procedure $P_2$ improves the power of the test and the test 
$T_{\mathcal{M}^3,P_2}$ has the greatest power. In this setting, one should 
  prefer the collection $\mathcal{M}^3$ to $\mathcal{M}^2$. This was 
  previously pointed out in Section \ref{ellipsoides} from a theoretical point 
  of view. Although $T_{\mathcal{M}^3,P_1}$ is conservative, it is a good 
  compromise for practical issues: it is very easy and fast to implement and 
  its performances are  good.

\section{Proofs of Theorem \ref{thrm_puissance}, Propositions 
  \ref{puissance_detection}, \ref{majoration_dependence}, 
  \ref{puissance_ellipsoide}, \ref{puissance_boules}, and \ref{majoration_minimax_ellipsoide}}\label{proof1} 
 
\begin{proof}[Proof of Theorem \ref{thrm_puissance}]

In a nutshell, we shall prove that conditionally to the design ${\bf X}$ the distribution of the test $T_{\alpha}$ is the same as the test introduced by Baraud \emph{et al.} \cite{baraud2003}. Hence, we may apply their non asymptotic upper bound for the power.~\\

\textbf{Distribution of $\phi_{m}({\bf Y},{\bf X})$ }.
First, we derive the distribution of the test statistic $\phi_{m}({\bf Y},{\bf X})$ under $\mathbb{P}_{\theta}$. 
The distribution of $Y$ conditionally to the set of 
variables $(X_{V\cup m})$ is of the form 
\begin{eqnarray} 
 Y = \sum_{i\in V\cup m} \theta_i^{V\cup m} X_i +   \epsilon^{V\cup m},  \label{decompositionyepsilon} 
\end{eqnarray} 
where the vector $\theta^{V\cup m}$ is constant and $\epsilon^{V\cup m}$ is a zero mean 
Gaussian variable independent of $X_{V\cup m}$, whose variance is $\var(Y|X_{V\cup m})$. As a consequence, $\|{\bf Y}-\Pi_{V\cup m}{\bf Y}\|_n^2$ is exactly 
$\|\Pi_{(V\cup m)^{\perp}}{\boldsymbol{\epsilon}}^{V\cup m}\|_n^2$, where 
$\Pi_{(V\cup m)^{\perp}}$ denotes the orthogonal projection along the 
space generated by $({\bf X}_i)_{i\in V\cup m}$. ~\\
Using the same decomposition of ${\bf Y}$ one simplifies the numerator of 
$\phi_m({\bf Y},{\bf X})$:
\begin{eqnarray} 
\left\|\Pi_{V\cup m}{\bf Y} - \Pi_{V}{\bf Y}\right\|_n^2 = \left\|\sum_{i\in 
  V\cup m}\theta_i^{V\cup m}({\bf X}_i-\Pi_V {\bf X}_i) + \Pi_{V^{\perp}\cap 
  (V\cup m)}{\boldsymbol{\epsilon}}^{V\cup m}\right\|_n^2, \nonumber 
\end{eqnarray} 
where $\Pi_{V^{\perp}\cap (V\cup m)}$ is the orthogonal projection onto the 
intersection between the space generated by $({\bf X}_i)_{i\in V\cup m}$ and 
the orthogonal of the space generated by $({\bf X}_i)_{i\in V}$. ~\\
For any $i \in m$, let us consider the conditional distribution of $X_i$ with respect to 
$ X_V$, 
\begin{eqnarray} 
 X_i = \sum_{j\in V} \theta_j^{V,i}  X_j +{\epsilon}_i^V. \label{decompositionxepsilon} 
\end{eqnarray} 
where $\theta_j^{V,i}$ are constants and $\epsilon_i^V$ is a 
 zero-mean normal gaussian random variable whose variance is $\var(X_i|X_V)$ 
 and which is independent of 
 $ X_V$. This enables us to express $${\bf X}_i-\Pi_{V}{\bf X}_i=\Pi_{V^{\perp}\cap 
 (V\cup m)}{\boldsymbol{\epsilon}}^V_i, \, \, \, \, \, \,\text{\bigskip for all }i\in m\ . $$ 
Therefore, we decompose $\phi_{m}({\bf Y}, {\bf X})$ in   
\begin{eqnarray}\label{decomposition_Tnm} 
\phi_{m}({\bf Y},{\bf X})= \frac{N_m \|\Pi_{V^{\perp} \cap (V\cup m)}\left(\sum_{i\in 
    m}\theta_i^{V\cup m}{\boldsymbol{\epsilon}}^{V}_i 
    +{\boldsymbol{\epsilon}}^{V\cup m}\right) \|_n^2}{D_m \|\Pi_{(V\cup 
    m)^{\perp}}{\boldsymbol{\epsilon}}^{V\cup m}\|_n^2}. 
\end{eqnarray} 
Let us define the random variable $Z_m^{(1)}$ and $Z_m^{(2)}$ where $Z_m^{(1)}$ refers to the 
numerator of (\ref{decomposition_Tnm}) divided by $N_m$ and $Z_m^{(2)}$ to the denominator 
divided by $D_m$. We now prove that $Z_m^{(1)}$ and 
$Z_m^{(2)}$ are independent.  ~\\~\\
The variables $({\boldsymbol{\epsilon}}_j^V)_{j\in m}$ are $\sigma\left({\bf X}_{V\cup m}\right)$-measurable
as linear combinations of elements in ${\bf X}_{V\cup m}$. Moreover, 
${\boldsymbol{\epsilon}}^{V\cup m}$ follows a zero mean normal distribution with covariance 
matrix $\var(Y|X_{V\cup m})I_n$ and is independent of ${\bf X}_{V\cup m}$. As a 
consequence, conditionally to ${\bf X}_{V\cup m}$, $Z_m^{(1)}$ and $Z_m^{(2)}$ are independent by 
Cochran's Theorem as they correspond to projections onto two sets orthogonal 
from each other.

As ${\boldsymbol{\epsilon}}_j^V$ is a 
  linear combination of the columns of ${\bf X}_{V\cup m}$, 
  $Z_m^{(1)}$ follows a non-central $\chi^2$ distribution conditionally to ${\bf X}_{V\cup m}$:
\begin{eqnarray} 
(Z_m^{(1)}|{\bf X}_{V\cup m}) \sim \var(Y|X_{V\cup m})\chi^2\left(\frac{\left\|\sum_{j\in 
      m}\theta_j^{V\cup m}\Pi_{(V\cup m) \cap V^{\perp} 
      }{\boldsymbol{\epsilon}}_j^{V} \right\|_n^2}{\var(Y|X_{V\cup m})},D_m\right).\nonumber 
\end{eqnarray} 
We denote $a_m^2({\bf X}_{V\cup m}):= \frac{\left\|\sum_{j\in 
      m}\theta_j^{V\cup m}\Pi_{(V\cup m) \cap V^{\perp} 
      }{\boldsymbol{\epsilon}}_j^{V} \right\|_n^2}{\var(Y|X_{V\cup m})}$ this non-centrality parameter. \\ \\

\textbf{Power of $T_{\alpha}$ conditionally to ${\bf X}_{V\cup m}$}. Conditionally to ${\bf X}_{V\cup m}$ our test statistic $\phi_{m}({\bf Y},{\bf X})$ is the same as that proposed by Baraud \emph{et al} \cite{baraud2003} with $n-d$ data and $\sigma^2=\var(Y|X_{V\cup m})$. Arguing as in their proof of Theorem 1, there exists some quantity $\bar{\Delta}_{m}(\delta)$ such that the procedure accepts the hypothesis with probability not larger than $\delta/2$ if $a_m^2({\bf X}_{V\cup m})>\bar{\Delta}_{m}(\delta)$:
\begin{eqnarray}\label{definition_delta_bar}
\bar{\Delta}_m(\delta) := 2.5\sqrt{1+K_m^2(U)}\sqrt{D_m\log\left(\frac{4}{\alpha_m\delta}\right)}\left(1+\sqrt{\frac{D_m}{N_m}}\right)+ \\2.5 \left[k_mK_m(U)\vee 5\right]\log\left(\frac{4}{\alpha_m\delta}\right)\left(1+\frac{2D_m}{N_m}\right) \nonumber ,
\end{eqnarray}
where $U_m:=\log(1/\alpha_m)$, $U:=\log(2/\delta)$, $k_m:=2\exp\left(4U_m/N_m\right)$, and 
$$K_m(u):= 1+2\sqrt{\frac{u}{N_m}}+2k_m\frac{u}{N_m}.$$
Consequently, we have 
\begin{eqnarray}\label{puissance1}
 \mathbb{P}_{\theta}\left(T_{\alpha}\leq 0|{\bf X}_{V\cup m}\right)\mathbf{1}\left\{a_{m}^2({\bf X}_{V\cup m})\geq \bar{\Delta}_m(\delta)\right\}\leq \delta/2.
\end{eqnarray}

Let derive the distribution of the non-central parameter $a_{m}({\bf X}_{V\cup m})$. First, we 
simplify the projection term as ${\boldsymbol{\epsilon}}_j^{V}$ is a linear 
combinations of elements of ${\bf X}_{V\cup m}$. 
$$\Pi_{(V\cup m) \cap V^{\perp} }{\boldsymbol{\epsilon}}_j^V = 
  \Pi_{V\cup m}{\boldsymbol{\epsilon}}_j^V - \Pi_V{\boldsymbol{\epsilon}}_j^V =\Pi_{V^{\perp}}{\boldsymbol{\epsilon}}_j^V.$$  Let us define $\kappa_m^2$ as 
\begin{eqnarray*} 
\kappa_m^2 := \frac{\var\left(\sum_{j\in m}\theta_j^{V\cup m}\epsilon_j^{V}\right)}{\var(Y|X_{V\cup m})}.  
\end{eqnarray*} 
As the variable $\sum_{j\in m}\theta_j^{V\cup m}{\boldsymbol{\epsilon}}_j^V$  is independent of 
${\bf X}_V$, and as almost surely the dimension of the vector space generated by 
${\bf X}_V$ is $d$,  we get
\begin{eqnarray*} 
\frac{\left\|\sum_{j\in m}\theta_j^{V\cup m}\Pi_{V^{\perp}} 
  {\boldsymbol{\epsilon}}_j^{V} \right\|_n^2}{\var(Y|X_{V\cup m})}\sim \kappa_m^2 \chi^2(n-d). 
\end{eqnarray*} 
Hence, applying for instance Lemma 1 in \cite{laurent98}, we get
\begin{eqnarray*}
\mathbb{P}_{\theta}\left[\frac{a^2_m({\bf X}_{V\cup m})}{\kappa_m^2}\geq (n-d)-2\sqrt{(n-d)U}\right]\leq \delta/2 .
\end{eqnarray*}

Let gather (\ref{puissance1}) with this last bound. If 
\begin{eqnarray}\label{definition_delta'}
\kappa_m^2\geq \Delta'_m(\delta):= \frac{\bar{\Delta}_m(\delta)}{(n-d)\left(1-2\sqrt{\frac{U}{n-d}}\right)}\ , 
\end{eqnarray}
then it holds that
\begin{eqnarray*} 
\mathbb{P}_{\theta}(T_{\alpha} \leq 0)& \leq&  \mathbb{P}_{\theta}\left(T_{\alpha}\leq 0,a^2_m({\bf X}_{V\cup m})>\bar{\Delta}_m(\delta)\right)+\mathbb{P}_{\theta}\left[a_m^2({\bf X}_{V\cup m})\leq \bar{\Delta}_m(\delta)\right]\\
& \leq & \mathbb{E}_{\theta}\left\{\mathbb{P}_{\theta}\left[T_{\alpha}\leq 0,a^2_m({\bf X}_{V\cup m})>\bar{\Delta}_m(\delta)|{\bf X}_{V\cup m}\right]\right\}+ \\ & & \mathbb{P}_{\theta}\left[\frac{a^2_m({\bf X}_{V\cup m})}{\kappa_m^2}\geq (n-d)-2\sqrt{(n-d)U}\right]\\
&\leq \delta.
\end{eqnarray*}

\textbf{Computation of $\kappa_m^2$}.
Let us now compute the quantity $\kappa_m^2$ in order to simplify Condition (\ref{definition_delta'}). Let first express $\var(Y|X_V)$ in terms of $\var(Y|X_{m\cup 
  V})$ using the decomposition (\ref{decompositionyepsilon}) of $Y$.  
\begin{eqnarray} 
\var(Y|X_V) &  = & \var\left(\sum_{j\in V\cup m}\theta_j^{V\cup m}X_j 
+\epsilon^{V\cup m}\left| X_V\right.\right) \nonumber \\ 
& = & \var\left(\sum_{j\in V\cup m}\theta_j^{V\cup m}X_j\left|X_V\right.\right) + 
\var\left(\epsilon^{V\cup m}\left|X_V\right.\right) \nonumber\\ 
& = & \var\left(\sum_{j\in V\cup m}\theta_j^{V\cup m}X_j\left|X_V\right.\right) + 
\var\left(Y\left|X_{V\cup m}\right.\right), 
\label{kappa1} 
\end{eqnarray} 
as ${\epsilon}^{V\cup m}$ is independent of ${X}_{V\cup m}$. Now using the definition of 
${\epsilon}_j^V$ in (\ref{decompositionxepsilon}), it turns out that 
\begin{eqnarray} 
\var\left(\sum_{j\in V\cup m}\theta_j^{V\cup m} X_j|X_V\right) & = & 
 \var\left(\sum_{j\in m}\theta_j^{V\cup m}X_j |X_V\right) \nonumber\\ 
 & = & \var\left(\sum_{j\in m}\theta_j^{V\cup m} \epsilon_j^V|X_V\right) \nonumber \\ 
 & = & \var\left(\sum_{j\in m}\theta_j^{V\cup m}\epsilon_j^V \right), \label{kappa2} 
\end{eqnarray}  
as the $(\epsilon_j^V)_{j\in m}$ are independent of $X_V$. Gathering formulae (\ref{kappa1}) and (\ref{kappa2}), we get 
\begin{eqnarray} 
\kappa_m^2 = \frac{\var(Y|X_V)-\var(Y|X_{V\cup m})}{\var(Y|X_{V\cup m})}. 
\end{eqnarray}~\\

Under Assumption $H_{\mathcal{M}}$, $U_m\leq N_m /10 $ for all 
$m\in\mathcal{M}$ and  
$U\leq N_m / 21$. Hence, the terms  $U/N_m$, $U_m/N_m$, $k_m$, and $K_m(U)$ behave like constants and it follows from (\ref{definition_delta'}) that $\Delta'(m)\leq \Delta(m)$, which concludes the proof.
 \end{proof}
 
\begin{proof}[Proof of Proposition \ref{puissance_detection}]
 
We first recall the classical upper bound for the binomial coefficient 
(see for instance (2.9) in \cite{massartflour}).
$$\log\left|\mathcal{M}(k,p)\right| = \log \left(_k^p\right) \leq k\log\left(\frac{ep}{k}\right).$$ 
As a consequence, $\log(1/\alpha_m)\leq \log(1/\alpha) + 
k\log\left(\frac{ep}{k}\right)$. Assumption (\ref{hypothese_detection}) with $L=21$
therefore implies Hypothesis $H_{\mathcal{M}}$. Hence, we are in position to apply the second result of Theorem 
\ref{thrm_puissance}. Moreover, the assumption on $n$ implies 
that $n\geq 21 k $ and $D_m/N_m$ is thus smaller than $1/20$ for any 
model $m$ in $\mathcal{M}(k,p)$. Formula (\ref{inegalite_puissance_hm}) in Theorem \ref{thrm_puissance} then 
translates into 
$$\triangle(m)\leq \frac{(1+\sqrt{0.05})L_1\left(\sqrt{k^2\log\left(\frac{ep}{k}\right)}+ 
 \sqrt{k\log\left(\frac{2}{\alpha\delta}\right)}\right)+1.1L_2\left(k\log\left(\frac{ep}{k}\right)+\log\left(\frac{2}{\alpha\delta}\right)\right)}{n},$$ 
and it follows that Proposition \ref{puissance_detection} holds. 
 
 \end{proof}
 
\begin{proof}[Proof of Proposition \ref{majoration_dependence}]
We fix the constant $L$ in Hypothesis (\ref{hypothese2majoration}) to be $21\log(4e)\vee C_2\log(4)$ where the universal constant $C_2$ is defined later in the proof. This choice of constants allows the procedure 
$\left[\sup_{1\leq i\leq p} \phi_{\{i\},\alpha/(2p)}\right]$ to satisfy Hypothesis $H_{\mathcal{M}}$. An argument similar to the proof of Proposition 
\ref{puissance_detection} allows to show easily that there exists a universal 
constant $C$ such that if we set 
\begin{eqnarray} 
\rho'^2_{1}& := & \frac{C\left(\log(p)+ 
  \log\left(\frac{4}{\alpha\delta}\right)\right)}{n} =  \frac{C}{n}\log\left(\frac{4p}{\alpha\delta}\right),\label{condition3_ulambda} 
\end{eqnarray} 
then  $\frac{\|\theta\|^2}{\var(Y)- 
 \|\theta\|^2}\geq \rho'^2_{1}$ implies that 
$\mathbb{P}_{\theta}\left(T_{\alpha}>0\right)\geq 1-\delta$. Here, the factor
 $4 $ in the logarithm comes from the fact that some weights $\alpha_m$ equal $\alpha/(2p)$. ~\\~\\
Let $\rho^2$  and $\lambda^2$ be two positive numbers such that 
$\frac{\lambda^2}{\var(Y)-\lambda^2}=\rho^2$ and let $\theta \in \Theta[1,p]$ 
such that $\|\theta\|^2=\lambda^2$. As $\corr(X_i,X_j)= c$ for any $i\neq j$, it follows that $\var(X_{p+1})=c+\frac{1-c}{p}$ and $\cov(Y,X_{p+1})^2=\|\theta\|^2\left[c+\frac{1-c}{p}\right]^2$.
$$\frac{\var(Y)-\var(Y|X_{p+1})}{\var(Y|X_{p+1})} = \frac{\left(c+(1-c)/p\right)\lambda^2}{\var(Y)-\left(c+(1-c)/p\right)\lambda^2}.$$    
We now apply Theorem \ref{thrm_puissance} to $\phi_{\{p+1\},\alpha/2}$ under $H_{\mathcal{M}}$. There exists a universal constant $C_2$ such that $\mathbb{P}_{\theta}\left(\phi_{\{p+1\},\alpha/2}>0\right)\geq 1-\delta$ if 
\begin{eqnarray*} 
\frac{\left(c+(1-c)/p\right)\lambda^2}{\var(Y)-\left(c+(1-c)/p\right)\lambda^2}\geq \frac{C_2}{n} \log\left(\frac{4}{\alpha\delta}\right). 
\end{eqnarray*}  
This last condition is implied by
\begin{eqnarray*}
 \frac{c\lambda^2}{\var(Y)-c\lambda^2}\geq \frac{C_2}{n} \log\left(\frac{4}{\alpha\delta}\right),
\end{eqnarray*}
which is equivalent to 
\begin{eqnarray} 
\frac{\lambda^2}{\var(Y)} \geq \frac{C_2}{cn + cC_2 \log\left(\frac{4}{\alpha\delta}\right) }\log\left(\frac{4}{\alpha\delta}\right).\label{condition_ulambda} 
\end{eqnarray} 
Let us assume that $c\geq \log\left(\frac{4}{\alpha\delta}\right)/\log\left(\frac{4p}{\alpha\delta}\right)$.  
As $n\geq 2C_2\log\left(\frac{4p}{\alpha\delta}\right)$ (Hypothesis 
(\ref{hypothese2majoration}) and definition of $L$), $nc\geq 
2C_2\log\left(\frac{4}{\alpha\delta}\right)$. As a consequence, Condition 
(\ref{condition_ulambda}) is implied by: 
\begin{eqnarray} 
\rho^2 \geq \frac{2C_2}{nc}\log\left(\frac{4}{\alpha\delta}\right).\label{condition2_ulambda} 
\end{eqnarray} 
Combining (\ref{condition3_ulambda}) and (\ref{condition2_ulambda}) allows  
to conclude that $\mathbb{P}_{\theta}\left(T_{\alpha}>0\right)\geq 1-\delta$ if  
\begin{eqnarray*} 
  \rho^2 \geq\frac{L}{n}\left(\log\left(\frac{4p}{\alpha\delta}\right)\bigwedge\frac{1}{c}\log\left(\frac{4}{\alpha\delta}\right)\right).  
\end{eqnarray*} 
\end{proof}
 
\begin{proof}[Proof of Proposition \ref{puissance_ellipsoide}] 
 We fix the constant $L$ to $42\log(80)$ in Hypothesis (\ref{condition_n_ellipsoide}). It follows that (\ref{condition_n_ellipsoide}) implies 
\begin{eqnarray}\label{condition_n2_ellipsoide}
 n\geq 42\left(\log\left(\frac{40}{\alpha}\right)\vee \log\left(\frac{2}{\delta}\right)\right).
\end{eqnarray}
First, we check that the test $T_{\alpha}$ satisfies Condition $H_{\mathcal{M}}$. 
As the dimension of each model is smaller than 
$n/2$, for any model $m$ in $\mathcal{M}$, $N_m$ is larger than 
$n/2$. Moreover, for any model $m$ in $\mathcal{M}$, $\alpha_m$ is larger than 
$\alpha/(2|\mathcal{M}|)$ and $|\mathcal{M}|$ is smaller than $n/2$. As a 
consequence, the first condition of $H_{\mathcal{M}}$ is implied by the 
inequality 
\begin{eqnarray}\label{condition1_hm_ellipsoide} 
n\geq 20\log\left(\frac{n}{\alpha}\right). 
\end{eqnarray} 
Hypothesis (\ref{condition_n2_ellipsoide}) implies that $n/2\geq 
20\log\left(\frac{40}{\alpha}\right)$. Besides, for any $n>0$ it holds that 
$n/2\geq 20\log\left(\frac{n}{40}\right)$. Combining these two lower bounds 
enables to obtain (\ref{condition1_hm_ellipsoide}). The second condition of 
$H_{\mathcal{M}}$ holds if $n\geq 42\log\left(\frac{2}{\delta}\right)$ which is 
a consequence of hypothesis (\ref{condition_n2_ellipsoide}). ~\\~\\
Let first consider the case $n<2p$ and let apply Theorem \ref{thrm_puissance} under Hypothesis 
$H_{\mathcal{M}}$ to $T_{\alpha}$. $\mathbb{P}_{\theta}\left(T_{\alpha}>0 \right)\geq 1-\delta$ 
for all $\theta\in \mathbb{R}^p$ such that 
\begin{equation}\label{condition1_ellipsoide} 
\exists i \in \{1,\ldots, [n/2]\}, \frac{\var(Y) - 
  \var\left(Y|X_{m_i}\right)}{\var\left(Y|X_{m_i}\right)}\geq C \frac{\sqrt{i\log\left(\frac{2[n/2]}{\alpha\delta}\right)}+ \log\left(\frac{2[n/2]}{\alpha\delta}\right)}{n},
\end{equation} 	
where $C$ is an universal constant. Let $\theta$ be an element of $\mathcal{E}_a(R)$ that satisfies
\begin{eqnarray*} 
\|\theta\|^2 \geq (1+C)\left(\var(Y|X_{m_i}) - \var(Y|X)\right) + (1+C)\var(Y|X) \frac{\sqrt{i\log\left(\frac{n}{\alpha\delta}\right)}+ \log\left(\frac{n}{\alpha\delta}\right)}{n}, 
\end{eqnarray*} 
for some $1\leq i\leq [n/2]$.
By Hypothesis (\ref{condition_n_ellipsoide}), it holds that
$$\frac{\sqrt{i\log\left(\frac{n}{\alpha\delta}\right)}+ 
  \log\left(\frac{n}{\alpha\delta}\right)}{n}\leq 1\ ,$$
for any $i$ between $1$ and $[n/2]$. It is then
straightforward to check that $\theta$
satisfies (\ref{condition1_ellipsoide}).~\\
As $\theta$ belongs to the set $\mathcal{E}_a(R)$,  
\begin{eqnarray*} 
\var(Y|X_{m_i}) - \var(Y|X) & = & a_{i+1}^2\var(Y|X) \sum_{j=i+1}^p 
\frac{\var(Y|X_{m_{j-1}})-\var(Y|X_{m_j})}{a_{i+1}^2\var(Y|X)} \\ 
& \leq & a_{i+1}^2\var(Y|X)R^2. 
\end{eqnarray*} 
Hence, if $\theta$ belongs to $\mathcal{E}_a(R)$ and satisfies  
\begin{eqnarray*} 
\|\theta\|^2 \geq (1+C) \var(Y|X)\left[\left(a_{i+1}^2 R^2 + 
  \frac{\sqrt{i\log\left(\frac{n}{\alpha\delta}\right)}}{n}\right) + 
  \frac{1}{n}\log\left(\frac{n}{\alpha\delta}\right) \right],
\end{eqnarray*} 
then $\mathbb{P}_{\theta}(T_{\alpha}\leq 0)\leq \delta$.
Gathering this condition for any $i$ between $1$ and $[n/2]$ allows to
conclude that if $\theta$ satisfies 
\begin{eqnarray*} 
\frac{\|\theta\|^2}{\var(Y)-\|\theta\|^2} \geq (1+C) \left[\inf_{1\leq i \leq [n/2]}\left(a_{i+1}^2 R^2 + 
  \frac{\sqrt{i\log\left(\frac{n}{\alpha\delta}\right)}}{n}\right) + 
  \frac{1}{n}\log\left(\frac{n}{\alpha\delta}\right) \right], 
\end{eqnarray*} 
then $\mathbb{P}_{\theta}(T_{\alpha}\leq 0)\leq \delta$.~\\~\\
Let us now turn to the case $n\geq 2p$. Let us consider $T_{\alpha}$
as the supremum of $p-1$ tests of level $\alpha/2(p-1)$ and one test of 
level $\alpha/2$. By 
considering the $p-1$ firsts tests, we obtain as in the previous case that 
$\mathbb{P}_{\theta}(T_{\alpha}\leq 0)\leq \delta$ if 
\begin{eqnarray*} 
\frac{\|\theta\|^2}{\var(Y)-\|\theta\|^2} \geq (1+C) \left[\inf_{1\leq i \leq (p-1)}\left(a_{i+1}^2R^2  + 
  \frac{\sqrt{i\log\left(\frac{p}{\alpha\delta}\right)}}{n}\right) + 
  \frac{1}{n}\log\left(\frac{p}{\alpha\delta}\right) \right]. 
\end{eqnarray*} 
On the other hand, using the last test statistic $\phi_{\mathcal{I},\alpha/2}$, 
$\mathbb{P}_{\theta}(T_{\alpha}\leq 0)\leq \delta$ if 
\begin{eqnarray*} 
\frac{\|\theta\|^2}{\var(Y)-\|\theta\|^2} \geq 
C\frac{\sqrt{p\log\left(\frac{2}{\alpha\delta}\right)}+ \log\left(\frac{2}{\alpha\delta}\right)}{n}. 
\end{eqnarray*} 
Gathering these two conditions allows to prove (\ref{majoration_ngrand}). 
 \end{proof}

\begin{proof}[Proof of Proposition \ref{puissance_boules}] 
 
The approach behind this proof is similar to the one for Proposition 
\ref{puissance_ellipsoide}. We fix the constant $L$ in Assumption \ref{condition_n_boules}, as in the previous proof.
Hence, the collection of models
$\mathcal{M}$ and the weights $\alpha_m$ satisfy hypothesis $H_{\mathcal{M}}$ as
in the previous proof.~\\~\\
Let us give a sharper upper bound on $|\mathcal{M}|$: 
\begin{eqnarray}\label{majoration_nbremodele} 
|\mathcal{M}|\leq 1+\log(n/2\wedge p)/\log(2)\leq \log(n\wedge 2p)/\log(2). 
\end{eqnarray} 
We deduce from (\ref{majoration_nbremodele}) that there exists a constant 
$L(\alpha,\delta)$ only depending on $\alpha$ and $\delta$ such that for all $m\in \mathcal{M}$, 
$$\log\left(\frac{1}{\alpha_m\delta}\right)\leq L(\alpha,\delta)\log\log 
(n\wedge p).$$ ~\\
 
First, let us consider the case $n<2p$. We apply Theorem \ref{thrm_puissance} under  
Assumption $H_{\mathcal{M}}$. As in the proof of Proposition 
\ref{puissance_ellipsoide}, we obtain that 
$\mathbb{P}_{\theta}(T_{\alpha}>0)\geq 1-\delta$ if 
\begin{eqnarray*} 
\frac{\|\theta\|^2}{\var(Y)-\|\theta\|^2} \geq 
L(\alpha,\delta)\left[\inf_{i\in \{2^j,j\geq 0\}\cap \{1,\ldots,[n/2]\}} \left(R^2(i+1)^{-2s} 
  +\frac{\sqrt{i\log\log n}}{n}\right) 
+\frac{\log\log n}{n}\right].
\end{eqnarray*} 
It is worth noting that $R^2i^{-2s}\leq \frac{\sqrt{i\log\log n}}{n}$ if and only if  
$$i\geq i^* = \left(\frac{R^2n}{\sqrt{\log\log n}}\right)^{2/(1+4s)}.$$  
Under the assumption on $R$, $i^*$ is larger than one. Let us distinguish 
between two cases. If there exists $i'$ in $\{2^j,j\geq 0\}\cap \{1,\ldots,[n/2]\}$ such that $i^*\leq 
i'$, one can take $i'\leq 2i^*$ and then 
\begin{eqnarray}\label{cas1_boule} 
\inf_{i\in\{2^j,j\geq 0\}\cap \{1,\ldots,[n/2]\}} \left(R^2i^{-2s} +\frac{\sqrt{i\log\log 
    n}}{n}\right) & \leq & 2\frac{\sqrt{i'\log\log n}}{n} \nonumber\\ 
& \leq & 2\sqrt{2}R^{2/(1+4s)}\left(\frac{\sqrt{\log \log n}}{n}\right)^{4s/(1+4s)}. 
\end{eqnarray} 
Else, we take $i'\in \{2^j,j\geq 0\}\cap \{1,\ldots,[n/2]\}$ such that $n/4\leq i'\leq 
n/2$. Since 
$i'\leq (i^*\wedge n/2)$ we obtain that 
\begin{eqnarray} \label{cas2_boule} 
\inf_{i\in  \{2^j,j\geq 0\}\cap \{1,\ldots,[n/2]\} } \left(R^2i^{-2s} +\frac{\sqrt{i\log\log 
    n}}{n}\right)  \leq 2R^2i'^{-2s}\leq 2R^2\left(\frac{n}{2}\right)^{-2s}. 
\end{eqnarray} 
Gathering inequalities (\ref{cas1_boule}) and (\ref{cas2_boule}) allows to 
 prove (\ref{majoration_npetit_boule}). ~\\~\\
We now turn to the case $n\geq 2p$. As in the proof of Proposition 
\ref{puissance_ellipsoide}, we divide the proof into two parts: first we give an 
upper bound of the power for the $|\mathcal{M}|-1$ first tests which define $T_{\alpha}$ and 
then we give an upper bound for the last test $\phi_{\mathcal{I},\alpha/2}$. Combining 
these two inequalities allows us to prove (\ref{majoration_ngrand_boule}). 
 \end{proof}
 
\begin{proof}[Proof of Proposition \ref{majoration_minimax_ellipsoide}]
 
We fix the constant $L$ in the assumption as in the two previous proofs. 
We first note that the assumption on $R^2$ implies that $D^*\geq 2$. 
As $N_m$ is larger than $n/2$,  the $\phi_{m_{D^*}}$ test clearly 
satisfies Condition $H_{\mathcal{M}}$. As a consequence, we may apply Theorem 
\ref{thrm_puissance}. Hence, $\mathbb{P}_{\theta}(T^*_{\alpha}\leq 0)\leq \delta$ for any 
$\theta$  such that 
\begin{eqnarray}\label{condition1_minimax_ellipsoide} 
  \frac{\var(Y)-\var(Y|X_{m_{D^*}})}{\var(Y|X_{m_{D^*}})}\geq 
  L(\alpha,\delta)\frac{\sqrt{D^*}}{n}.
\end{eqnarray} 
 Now, we use the same sketch as in the proof of Proposition 
\ref{puissance_ellipsoide}. For any $\theta\in \mathcal{E}_a(R)$, Condition 
(\ref{condition1_minimax_ellipsoide}) is equivalent to: 
\begin{equation} 
\|\theta\|^2\geq \left(\var(Y|X_{m_{D^*}})- \var(Y|X) \right)\left(1+L(\alpha,\delta)\frac{\sqrt{D^*}}{n}\right)+\var(Y|X)L(\alpha,\delta)\frac{\sqrt{D^*}}{n}. \label{condition2_minimax_ellipsoide}
\end{equation} 
Moreover, as $\theta$ belongs to $\mathcal{E}_a(R)$, $$\var(Y|X_{m_{D^*}})- 
\var(Y|X) \leq a_{D^*+1}^2R^2\var(Y|X)\leq a_{D^*}^2\var(Y|X)R^2.$$ 
As $\sqrt{D^*}/n$ is smaller than one, Condition
(\ref{condition2_minimax_ellipsoide}) is implied by 
\begin{eqnarray*} 
\frac{\|\theta\|^2}{\var(Y)-\|\theta\|^2}\geq (1+L(\alpha,\delta))\left(a_{D^*}^2R^2+\frac{\sqrt{D^*}}{n}\right). 
\end{eqnarray*} 
As $a_{D^*}^2R^2$ is smaller than $\frac{\sqrt{D^*}}{n}$ which is smaller 
$\sup_{1\leq i\leq p}\left[\frac{\sqrt{i}}{n}\wedge a_i^2R^2\right]$, it turns 
out that $\mathbb{P}_{\theta}(T^*_{\alpha}=0)\leq \delta$ for any $\theta$ 
belonging to $\mathcal{E}_a(R)$ such that 
\begin{eqnarray*} 
\frac{\|\theta\|^2}{\var(Y)-\|\theta\|^2}\geq 2(1+L(\alpha,\delta))\sup_{1\leq i\leq p}\left[\frac{\sqrt{i}}{n}\wedge a_i^2R^2\right]. 
\end{eqnarray*} 
\end{proof} 

\section{Proofs of Theorem \ref{minoration_arete}, Propositions  
  \ref{minoration_puissance_1test}, \ref{minoration_1test}, \ref{minoration_dependence}, 
  \ref{spatial}, \ref{vitesse_minimax_ellipsoide}, \ref{minimax_lower_nested_linear}, and \ref{minoration_minimax_collection_ellipsoide}}\label{proof2} 

Throughout this section, we shall use the notations $\eta := 2(1-\alpha - \delta)$ and 
$\mathcal{L}(\eta):=  \frac{\log(1+2\eta^2)}{2}$. 
 
\begin{proof}[Proof of Theorem \ref{minoration_arete}]

This proof follows the general method for obtaining lower bounds described 
in Section 7.1 in Baraud \cite{baraudminimax}. We first remind the reader of the 
main arguments of the approach applied to our model. Let $\rho$ be some positive 
number and $\mu_{\rho}$ be some probability measure on 
\begin{eqnarray*} 
\Theta[k,p,\rho] := \left\{\theta\in  \Theta[k,p], \frac{\|\theta \|^2}{\var(Y)-\|\theta 
\|^2}= \rho\right\}. 
\end{eqnarray*} 
We define $\mathbb{P}_{\mu_{\rho}}=\int \mathbb{P}_\theta \text{d}\mu_{\rho}(\theta)$ and $\Phi_{\alpha}$ 
the set of level-$\alpha$ tests of the hypothesis "$\theta=0$". Then, 
\begin{eqnarray} 
  \beta_I(\Theta[k,p,\rho]) & \geq & \inf_{\phi_\alpha \in \Phi_\alpha} 
\mathbb{P}_{\mu_\rho}[\phi_{\alpha}=0]\nonumber \\  
& \geq & 1-\alpha - \sup_{A,\ \mathbb{P}_0(A)\leq \alpha}|\mathbb{P}_{\mu_\rho}(A)-\mathbb{P}_{0}(A)| 
\nonumber\\ 
& \geq & 1-\alpha-\frac{1}{2}\|\mathbb{P}_{\mu_{\rho}}-\mathbb{P}_0 \|_{TV},\label{minoration_puissance} 
\end{eqnarray} 
where $\|\mathbb{P}_{\mu_\rho}- \mathbb{P}_0\|_{TV}$ denotes the total variation norm between the 
probabilities $\mathbb{P}_{\mu_{\rho}}$ and $\mathbb{P}_0$. If we suppose that $\mathbb{P}_{\mu_{\rho}}$ 
is absolutely continuous with respect to $\mathbb{P}_0$, we can upper bound the 
norm in total variation between these two probabilities as follows. We define  
$$L_{\mu_{\rho}}({\bf Y},{\bf X}) := \frac{\text{d}\mathbb{P}_{\mu_{\rho}}}{\text{d}\mathbb{P}_0}(\bf{Y},\bf{X}).$$ 
Then, we get the upper bound  
\begin{eqnarray*} 
\|\mathbb{P}_{\mu_{\rho}}-\mathbb{P}_0\|_{TV} & = & \int |L_{\mu_{\rho}}({\bf Y},{\bf X})- 1 
|\text{d}\mathbb{P}_0({\bf Y},{\bf X})\\ 
& \leq & \left(\mathbb{E}_0\left[L^2_{\mu_{\rho}}({\bf Y},{\bf X})\right]-1\right)^{1/2}. 
\end{eqnarray*} 
Thus, we deduce from (\ref{minoration_puissance}) that 
\begin{eqnarray*} 
\beta_I(\Theta[k,p,\rho]) \geq 1 - \alpha - \frac{1}{2}\left(\mathbb{E}_0\left[L^2_{\mu_{\rho}}({\bf Y},{\bf X})\right]-1\right)^{1/2}. 
\end{eqnarray*} 
If we find a number $\rho^*=\rho^*(\eta)$ such that 
\begin{eqnarray}\label{argument_principal_minimax} 
\log\left(\mathbb{E}_0\left[L^2_{\mu_{\rho^*}}({\bf Y},{\bf X})\right]\right)\leq \mathcal{L}(\eta), 
\end{eqnarray} 
then for any $\rho\leq \rho^*$, 
\begin{eqnarray*} 
\beta_I(\Theta[k,p,\rho]) \geq 1 - \alpha - \frac{\eta}{2} = \delta. 
\end{eqnarray*}

To apply this method, we first have to define a suitable prior 
$\mu_{\rho}$ on $\Theta[k,p,\rho]$. Let $\widehat{m}$ be some random variable uniformly 
distributed over $\mathcal{M}(k,p)$ and for each $m\in \mathcal{M}(k,p)$, let 
$\epsilon^{m}=(\epsilon^m_j)_{j\in m}$ be a sequence of independent Rademacher random 
variables. We assume that for all $m\in \mathcal{M}(k,p)$, $\epsilon^m$ and 
$\widehat{m}$ are independent. Let $\rho$ be given and $\mu_{\rho}$ be the 
distribution of the random variable 
$\widehat{\theta}= \sum_{j\in\widehat{m}} \lambda\epsilon_j^{\widehat{m}}e_j$ where 
\begin{eqnarray*} 
  \lambda^2 := \frac{\var(Y) \rho^2}{k(1+\rho^2)}, 
\end{eqnarray*} 
and where $(e_j)_{j\in \mathcal{I}}$ is the orthonormal family of vectors of 
$\mathbb{R}^{p}$ defined by 
$$ (e_j)_i =1 \text{ if } i=j \text{ and } (e_i)_j=0 \text{ otherwise}. $$ 
 
Straightforwardly, $\mu_{\rho}$ is supported by 
$\Theta[k,p,\rho]$.  
For any $m$ in $\mathcal{M}(k,p)$ and any vector $(\zeta^m_j)_{j\in m}$ with 
values in $\{-1;1\}$, let  
$\mu_{m,\zeta^m,\rho}$ be the Dirac measure on $\sum_{j\in m} 
\lambda\zeta_j^{m}e_j$. For any $m$ in $\mathcal{M}(k,p)$, 
$\mu_{m,\rho}$ denotes the distribution of the random variable 
$\sum_{j\in m} \lambda\zeta_j^{m}e_j$ where $(\zeta_j^m)$ is a 
sequence of independent Rademacher random variables. These definitions easily imply  
$$L_{\mu_{\rho}}({\bf Y},{\bf X}) = 
\frac{1}{\left(_k^p\right)}\sum_{m\in\mathcal{M}(k,p)}L_{\mu_{m,\rho}}({\bf Y},{\bf X})= 
\frac{1}{2^{k}\left(_k^p\right)}\sum_{m\in\mathcal{M}(k,p)}\sum_{\zeta^m \in \{-1,1\}^k}L_{\mu_{m,\zeta^m\rho}}({\bf Y},{\bf X}).$$ 
We aim at bounding the quantity $\mathbb{E}_0(L_{\mu_\rho}^2)$ and obtaining an inequality of the form (\ref{argument_principal_minimax}). First, we work out $L_{\mu_{m,\zeta^m,\rho}}$:  
 
\begin{eqnarray} 
L_{\mu_{m,\zeta^m,\rho}}({\bf Y},{\bf X}) & = & 
\left[\left(\frac{1}{1-\frac{\lambda^2k}{\var(Y)}}\right)^{n/2}\exp\left(-\frac{\|{\bf Y}\|_n^2}{2}\frac{\lambda^2k}{\var(Y)(\var(Y)-\lambda^2k)}\right.\right.\nonumber\\ 
& + & \left.\left. \lambda\sum_{j\in m}\zeta^m_j 
  \frac{<{\bf Y},{\bf X}_j>_n}{\var(Y)-\lambda^2k} - \lambda^{2}\sum_{j,j'\in m}\zeta^m_j\zeta^m_{j'}\frac{<{\bf X}_j,{\bf X}_{j'}>_n}{2(\var(Y)-\lambda^2k)} \right)\right],  
\end{eqnarray} 
where $<.>_n$ refers to the canonical inner product in $\mathbb{R}^n$. 
 
Let us fix $m_1$ and $m_2$ in $\mathcal{M}(k,p)$ and two vectors $\zeta^1$ and 
$\zeta^2$ respectively associated to $m_1$ and $m_2$. We aim at computing 
the quantity 
$\mathbb{E}_0\left(L_{\mu_{m_1,\zeta^1,\rho}}({\bf Y},{\bf X})L_{\mu_{m_2,\zeta^2,\rho}}({\bf Y},{\bf X})\right)$. 
First, we decompose the set $m_1\cup m_2$ into 
four sets (which possibly are empty): $m_1\setminus m_2$, $m_2\setminus m_1$, $m_3$, and $m_4$, where $m_3$ and $m_4$ are defined by: 
\begin{eqnarray} 
m_3 & := & \left\{j\in m_1\cap m_2|\zeta^1_j =\zeta^2_j \right\}\nonumber\\ 
m_4 & := & \left\{j\in m_1\cap m_2|\zeta^1_j =-\zeta^2_j \right\}.\nonumber 
\end{eqnarray} 
 
For the sake of simplicity, we reorder the elements of 
$m_1\cup m_2$ from $1$ to $|m_1\cup m_2|$ such that the first elements 
belong to $m_1 \setminus m_2$, then to $m_2\setminus 
m_1$ and so on. Moreover, we define the vector $\zeta 
\in\mathbb{R}^{|m_1\cup m_2|}$ such that $\zeta_j=\zeta^1_j$ if $j\in 
m_1$ and $\zeta_j=\zeta^2_j$ if $j\in m_2 \setminus  m_1$. Using these 
notations, we compute the expectation of $L_{m_1,\zeta^1,\rho}({\bf Y},{\bf X})L_{m_2,\zeta^2,\rho}({\bf Y},{\bf X})$. 
 
\begin{eqnarray} \label{esperance_m_xi_rho}
\mathbb{E}_{0}\left(L_{\mu_{m_1,\zeta^1,\rho}}({\bf Y},{\bf X})L_{\mu_{m_2,\zeta^2,\rho}}({\bf Y},{\bf X})\right)= 
\left(\frac{1}{\var(Y)(1-\frac{\lambda^2k}{\var(Y)})^2}\right)^{n/2}|A|^{-n/2}\ ,
\end{eqnarray} 
where $|.|$ refers to the determinant and $A$ is a symmetric square matrix of size $|m_1\cup m_2|+1$ such that: 
\begin{eqnarray*} 
A[1,j] :=\left\{\begin{array}{ccc} 
 \frac{\var(Y)+\lambda^2k}{\var(Y)(\var(Y)-\lambda^2k)} & \text{if} & j=1\\ 
 -\frac{\lambda\zeta_{j-1}}{\var(Y)-\lambda^2k} & \text{if} & (j-1)\in 
 m_1\triangle m_2 \\ 
 -2\frac{\lambda\zeta_{j-1}}{\var(Y)-\lambda^2k} & \text{if} & (j-1)\in 
 m_3 \\ 
 0 & \text{if} & (j-1)\in 
 m_4\ , \\ 
\end{array} \right. 
\end{eqnarray*} 
where $m_1\triangle m_2$ refers to $(m_1\cup m_2)\setminus (m_1\cap m_2)$. For any $i>1$ and $j>1$, $A$ satisfies 
\begin{eqnarray*} 
A[i,j] := \left\{\begin{array}{ccc} 
\lambda^2 \frac{\zeta_{i-1}\zeta_{j-1}}{\var(Y) -\lambda^2k} + \delta_{i,j} & \text{if} & (i-1,j-1)\in (m_1\setminus m_2)\times m_1\\
\lambda^2 \frac{\zeta_{i-1}\zeta_{j-1}}{\var(Y) -\lambda^2k} + \delta_{i,j} & \text{if} & (i-1,j-1)\in
(m_2\setminus m_1)\times (m_2\setminus m_1\cup m_3) \\ 
-\lambda^2 \frac{\zeta_{i-1}\zeta_{j-1}}{\var(Y) -\lambda^2k} & \text{if} & (i-1,j-1)\in
(m_2\setminus m_1)\times m_4 \\ 
2\lambda^2 \frac{\zeta_{i-1}\zeta_{j-1}}{\var(Y) -\lambda^2k}+ \delta_{i,j} & \text{if} & (i-1,j-1)\in \left[m_3\times m_3\right]\cup \left[m_4\times m_4 \right]\\ 
0 & \text{else}, & 
\end{array}\right.  \ ,
\end{eqnarray*} 
where $\delta_{i,j}$ is the indicator function of $i=j$.

After some linear transformation on the lines of the matrix $A$, it is 
possible to express its determinant into 
$$|A| = \frac{\var(Y)+\lambda^2k}{\var(Y)(\var(Y)-\lambda^2k)}\left|I_{|m_1\cup m_2|}+C\right|,$$ 
where $I_{|m_1\cup m_2|}$ is the identity matrix of size $|m_1\cup m_2|$. 
$C$ is a symmetric matrix of size $|m_1\cup m_2|$ such that for any $(i,j)$,  
\begin{eqnarray*} 
C[i,j] = \zeta_i\zeta_jD[i,j] 
\end{eqnarray*} 
and $D$ is a block symmetric matrix defined by 
\begin{eqnarray*} 
D :=  
\left[\begin{array}{cccc}  
\frac{\lambda^4k}{\var^2(Y)-\lambda^4k^2} & 
\frac{-\lambda^2\var(Y)}{\var^2(Y)-\lambda^4k^2} & 
\frac{-\lambda^2}{\var(Y)+\lambda^2k} & 
\frac{\lambda^2}{\var(Y)-\lambda^2k}  
\\ \frac{-\lambda^2\var(Y)}{\var^2(Y)-\lambda^4k^2} & 
\frac{\lambda^4k}{\var^2(Y)-\lambda^4k^2} & \frac{-\lambda^2}{\var(Y)+\lambda^2k} & 
\frac{-\lambda^2}{\var(Y)-\lambda^2k}  
\\ \frac{-\lambda^2}{\var(Y)+\lambda^2k} 
&\frac{-\lambda^2}{\var(Y)+\lambda^2k} 
&\frac{-2\lambda^2}{\var(Y)+\lambda^2k} & 0  
\\ \frac{\lambda^2}{\var(Y)-\lambda^2k} & 
\frac{-\lambda^2}{\var(Y)-\lambda^2k} & 0 & \frac{2\lambda^2}{\var(Y)-\lambda^2k} 
\end{array}\right] \, .
\end{eqnarray*} 
Each block corresponds to one of the four previously defined subsets of $m_1\cup m_2$ 
(i.e. $m_1 \setminus m_2$, $m_2 \setminus m_1$, $m_3$, and $m_4$). The matrix $D$ 
is of rank at most four. By computing its non-zero eigenvalues, it is then 
straightforward to derive the determinant of $A$ 
\begin{eqnarray} 
  |A| = \frac{\left[\var(Y)-\lambda^2(2|m_3|-|m_1\cap m_2|)\right]^2}{\var(Y)(\var(Y)-\lambda^2k)^2}.  \nonumber 
\end{eqnarray} 
Gathering this equality with (\ref{esperance_m_xi_rho}) yields 
\begin{eqnarray}\label{expression_vraisemblance} 
\mathbb{E}_{0}\left(L_{\mu_{m_1,\zeta^1,\rho}}({\bf Y},{\bf X})L_{\mu_{m_2,\zeta^2,\rho}}({\bf Y},{\bf X})\right) = \left[\frac{1}{1-\frac{\lambda^2(2|m_3|-|m_1\cap m_2|)}{\var(Y) }}\right]^n. 
\end{eqnarray} 
 
Then, we take the expectation with respect to $\zeta^1$, 
$\zeta^2$, $m_1$ and $m_2$. When $m_1$ and $m_2$ are fixed the expression (\ref{expression_vraisemblance}) depends on $\zeta^1$ and $\zeta^2$ only towards the cardinality of $m_3$. As $\zeta^1$ and $\zeta^2$ correspond to independent Rademacher variables, the random variable $2|m_3|-|m_1\cap m_2|$ follows the distribution of $Z$, a sum of $|m_1\cap m_2|$ independent Rademacher variables and 
\begin{eqnarray} 
\mathbb{E}_0(L_{\mu_{m_1,\rho}}({\bf Y},{\bf X})L_{\mu_{m_2,\rho}}({\bf Y},{\bf X})) & = & \mathbb{E}_0\left[\frac{1}{1-\frac{\lambda^2 Z}{\var(Y) }}\right]^n \label{produit_vraisemblance}. 
\end{eqnarray} 
 
When $Z$ is non-positive, this expression is smaller than one. Alternatively, when $Z$ is non-negative:  
\begin{eqnarray*} 
\left[\frac{1}{1-\frac{\lambda^2 Z}{\var(Y) }}\right]^n & = & 
\exp\left(n\log\left(\frac{1}{1-\frac{\lambda^2Z}{\var(Y) }}\right)\right) \nonumber\\ 
& \leq & \exp\left[n\frac{\frac{\lambda^2Z}{\var(Y) }}{1 - 
    \frac{\lambda^2Z}{\var(Y) }}\right] \nonumber \\ 
& \leq & \exp\left[n\frac{\frac{\lambda^2Z}{\var(Y) }}{1 - 
    \frac{\lambda^2k}{\var(Y)}}\right] , 
\end{eqnarray*} 
as $\log(1+x)\leq x$ and as $Z$ is smaller than $k$. We define an event $\mathbb{A}$ 
such that $\{Z>0\} \subset \mathbb{A} \subset \{Z\geq 0 \}$ and $\mathbb{P}(\mathbb{A}) = 
\frac{1}{2}$.  This is always possible as the random variable $Z$ is 
symmetric. As a consequence, on the event $\mathbb{A}^c$, the quantity 
(\ref{produit_vraisemblance}) is smaller or equal to one. All in all, we bound (\ref{produit_vraisemblance}) by: 
\begin{eqnarray}\label{majoration_indicatrice_globale} 
\mathbb{E}_0(L_{\mu_{m_1,\rho}}({\bf Y},{\bf X})L_{\mu_{m_2,\rho}}({\bf Y},{\bf X})) \leq \frac{1}{2} 
+\mathbb{E}_0\left[\mathbf{1}_{\mathbb{A}} \exp\left[n\frac{\frac{\lambda^2Z}{\var(Y) }}{1 - 
    \frac{\lambda^2k}{\var(Y)}}\right]  \right], 
\end{eqnarray} 
where $\mathbf{1}_{\mathbb{A}}$ is the indicator function of the event $\mathbb{A}$. We now apply H\"{o}lder's inequality with a parameter 
$v \in ]0;1]$, which will be fixed later. 
 
\begin{eqnarray}\label{majoration_indicatrice} 
\mathbb{E}_0\left[\textbf{1}_{\mathbb{A}} \exp\left[n\frac{\frac{\lambda^2Z}{\var(Y) }}{1 - 
    \frac{\lambda^2k}{\var(Y)}}\right]  \right] & \leq & 
    \mathbb{P}(\mathbb{A})^{1-v} \left[\mathbb{E}_0 \exp\left(\frac{n}{v}\frac{\frac{\lambda^2Z}{\var(Y) }}{1 - 
    \frac{\lambda^2k}{\var(Y)}}\right)\right]^v\nonumber\\ 
& \leq & \left(\frac{1}{2}\right)^{1-v}\left[\cosh\left(\frac{n\lambda^2}{v 
    (\var(Y) - \lambda^2k)} \right)\right]^{|m_1\cap m_2|v}. 
\end{eqnarray} 
Gathering inequalities (\ref{majoration_indicatrice_globale}) and 
(\ref{majoration_indicatrice}) yields 	 
\begin{eqnarray*} 
\mathbb{E}_0\left[L^2_{\mu_{\rho}}({\bf Y},{\bf X}) \right] \leq \frac{1}{2} + 
    \left(\frac{1}{2}\right)^{1-v}\frac{1}{\left(_k^p\right)^2}\sum_{m_1,m_2\in \mathcal{M}(k,p)}\cosh\left(\frac{n\lambda^2}{v 
    (\var(Y) - \lambda^2k)} \right)^{|m_1\cap m_2|v}. 
\end{eqnarray*} 
 
Following the approach of Baraud \cite{baraudminimax} in Section $7.2$, we note that 
if $m_1$ and $m_2$ are taken uniformly and independently in $\mathcal{M}(k,p)$, 
then $|m_1\cap m_2|$ is distributed as a Hypergeometric distribution with 
parameters $p$, $k$, and $k/p$. Thus, we derive that  
\begin{eqnarray} 
\mathbb{E}_0\left[L^2_{\mu_{\rho}}({\bf Y},{\bf X}) \right] \leq \frac{1}{2} + 
    \left(\frac{1}{2}\right)^{1-v}\mathbb{E}\left(\cosh\left(\frac{n\lambda^2}{v 
    (\var(Y) - \lambda^2k)} \right)^{vT}\right)\label{majoration_lmurho} 
\end{eqnarray} 
where $T$ is a random variable distributed according to a Hypergeometric 
distribution with parameters $p$, $k$ and $k/p$. We know from 
Aldous (p.173) \cite{aldous85} that $T$ has the same distribution as 
the random variable $\mathbb{E}(W|\mathcal{B}_p)$ where $W$ is binomial random 
variable of parameters $k$, $k/p$ and $\mathcal{B}_p$ some suitable 
$\sigma$-algebra. By a convexity argument, we then upper bound (\ref{majoration_lmurho}).

\begin{eqnarray*} 
\mathbb{E}_0\left[L^2_{\mu_{\rho}}({\bf Y},{\bf X}) \right]& \leq & \frac{1}{2} + 
    \left(\frac{1}{2}\right)^{1-v}\mathbb{E}\left(\cosh\left(\frac{n\lambda^2}{v 
    (\var(Y) - \lambda^2k)} \right)^{vW}\right) \nonumber\\ 
& = & \frac{1}{2} + 
\left(\frac{1}{2}\right)^{1-v}\left(1+\frac{k}{p}\left(\cosh\left(\frac{n\lambda^2}{v 
   (\var(Y) -\lambda^2k)}\right)^{v}-1\right) \right)^k\nonumber\\ 
& = & \frac{1}{2} + 
\left(\frac{1}{2}\right)^{1-v}\exp\left[k\log\left(1+\frac{k}{p}\left(\cosh\left(\frac{n\lambda^2}{v 
   (\var(Y) -\lambda^2k)}\right)^{v}-1\right)\right) \right]. 
\end{eqnarray*} 
 
To get the upper bound on the total variation distance appearing in (\ref{minoration_puissance}), we aim at constraining 
this last expression to be smaller than $1+\eta^2$. This is equivalent to the 
following inequality: 
\begin{eqnarray}\label{ineq_eta} 
2^{v} \exp\left[k\log\left(1+\frac{k}{p}\left(\cosh\left(\frac{n\lambda^2k}{v 
    k (\var(Y) -\lambda^2k)}\right)^{v}-1\right)\right) \right]\leq 1+2\eta^2\ . 
\end{eqnarray}
We now choose $v = \frac{\mathcal{L}(\eta)}{\log(2)}\wedge 1$. If $v$ is 
strictly smaller than one, then (\ref{ineq_eta}) is equivalent to: 
 
\begin{eqnarray}\label{ineq_eta2} 
k\log\left[1+\frac{k}{p}\left(\cosh\left(\frac{n\lambda^2k}{v 
  k (\var(Y) -\lambda^2k)}\right)^{v}-1\right)\right] \leq \frac{\log(1+2\eta^2)}{2}. 
\end{eqnarray} 
It is straightforward to show that this 
last inequality also implies (\ref{ineq_eta}) if $v$ equals one. We now suppose that 
\begin{eqnarray}\label{condition_complexe_minoration} 
\frac{n\lambda^2}{v  (\var(Y) -\lambda^2k)}\leq 
\log\left((1+u)^{\frac{1}{v}}+\sqrt{(1+u)^{\frac{2}{v}}-1}\right), 
\end{eqnarray} 
where $u=\frac{p \mathcal{L}(\eta)}{k^2}$. Using the classical equality 
$\cosh\left[\log(1+x+\sqrt{2x+x^2})\right] =1+x$ with $x= (1+u)^{\frac{1}{v}} -1$, we 
deduce that inequality (\ref{condition_complexe_minoration}) implies 
(\ref{ineq_eta2}) because
\begin{eqnarray} 
k\log\left(1+\frac{k}{p}\left(\cosh\left(\frac{n\lambda^2k}{v 
  k (\var(Y) -\lambda^2k)}\right)^{v}-1\right)\right) & \leq & 
  k\log\left(1+\frac{k}{p}u\right) \nonumber\\ 
& \leq & \frac{k^2}{p}u \leq \mathcal{L}(\eta)\ .\nonumber  
\end{eqnarray} 
For any $\beta\geq 1$ and any $x>0$, it holds that $(1+x)^{\beta}\geq 1 +\beta 
x$. As $\frac{1}{v}\geq 1$, Condition (\ref{condition_complexe_minoration}) is 
implied by: 
$$\frac{\lambda^2k}{\var(Y)-\lambda^2k}\leq \frac{kv}{n} \log\left(1 +\frac{u}{v}+\sqrt{\frac{2u}{v}}\right).$$ 
One then combines the previous inequality with the definitions of $u$ and 
$v$ to obtain the upper bound
\begin{eqnarray*} 
\frac{\lambda^2k}{\var(Y)-\lambda^2k} \leq \frac{k}{n}\left(\frac{\mathcal{L}(\eta)}{\log(2)}\wedge 1\right) 
\log\left(1 +\frac{p(\log(2)\vee 
  \mathcal{L}(\eta))}{k^2}+\sqrt{\frac{2p(\log(2)\vee \mathcal{L}(\eta))}{k^2}}\right)\ . 
\end{eqnarray*} 
For any $x$ positive and any $u$ between 0 and 
1, $\log(1+ux)\geq u\log(1+x)$. As a consequence, the previous inequality is implied by: 
\begin{eqnarray*} 
 \frac{\lambda^2k}{\var(Y)-\lambda^2k} & \leq & 
 \frac{k}{n}\left(\frac{\mathcal{L}(\eta)}{\log(2)}\wedge 
 1\right)\left(\left[\mathcal{L}(\eta)\vee \log(2)\right]\wedge1\right) \log\left(1 
 +\frac{p}{k^2}+\sqrt{\frac{2p}{k^2}}\right) \nonumber\\ 
& = & \frac{k}{n}\left(\mathcal{L}(\eta)\wedge 1\right) \log\left(1 
 +\frac{p}{k^2}+\sqrt{\frac{2p}{k^2}}\right)\ . 
\end{eqnarray*} 
To resume, if we take $\rho^2$ smaller than (\ref{enonce_minoration}), then  
$$\beta_I\left(\Theta[k,p,\rho]\right)\geq \delta\ .$$ 
Besides, the lower bound is strict if $\rho^2$ is strictly smaller than (\ref{enonce_minoration}). 
To prove the second part of the theorem, one has to observe that $\alpha + 
\delta\leq 53\% $ implies that $\mathcal{L}(\eta) \geq \frac{1}{2}$.  
 \end{proof}

\begin{proof}[Proof of Proposition \ref{minoration_1test}]
Let us first assume that the covariance matrix of $X$ is the identity. We argue as in the proof of Theorem \ref{minoration_arete} taking $k=p$. The sketch of the proof remains unchanged except that we slightly modify the last 
part. Inequality (\ref{ineq_eta2}) becomes 
$$pv\log\left(\cosh\left(\frac{n\lambda^2p}{vp(\var(Y)-\lambda^2p)}\right)\right)\leq \mathcal{L}(\eta),$$ 
where we recall that $v = \frac{\mathcal{L}(\eta)}{\log 2}\wedge 1$. For all 
$x\in \mathbb{R}$, $\cosh(x)\leq \exp(x^2/2)$. 
Consequently, the previous inequality is implied by 
\begin{eqnarray*} 
\frac{\lambda^2p}{\var(Y)-\lambda^2p}\leq \sqrt{2v\mathcal{L}(\eta)}\frac{\sqrt{p}}{n}, 
\end{eqnarray*} 
and the result follows easily. 

If we no longer assume that the covariance matrix $\Sigma$ is the identity, we orthogonalize the sequence $X_i$ thanks to Gram-Schmidt process. Applying the previous argument to this new sequence of covariates allows to conclude.
 \end{proof}
 
\begin{proof}[Proof of Proposition \ref{minoration_puissance_1test}] 
 Let define the constant $L(\alpha,\delta)$ involved in the condition:
\begin{eqnarray*}
 L(\alpha,\delta):=\sqrt{\log\left(1+8(1-\alpha-\delta)^2\right)}\left[1\wedge\sqrt{\log\left(1+8(1-\alpha-\delta)^2\right)/(2\log 2)}\right]
\end{eqnarray*}

Let us apply proposition \ref{minoration_1test}. For any $\rho\leq L(\alpha,\delta)\frac{\sqrt{D_m}}{{n}}$ and any $\varsigma>0$ there exists 
some $\theta\in S_m$ such that 
$\frac{\|\theta\|^2}{\var(Y)-\|\theta\|^2}=\rho^2$ and 
$\mathbb{P}_{\theta}(\phi_{m,\alpha}\leq 0) \geq \delta-\varsigma$.  In the proof of Theorem 
\ref{thrm_puissance}, we have shown in (\ref{decomposition_Tnm}) and following 
equalities that the distribution of the test statistic $\phi_m$ only 
depends on the quantity $\kappa_m^2 =\frac{\var(Y)-\var(Y|X_m)}{\var(Y|X_m)}$. Let $\theta'$ 
be an element of $S_m$ such that $\kappa_m^2= \rho^2$. The 
distribution of $\phi_m$ under $\mathbb{P}_{\theta'}$ is the same as its 
distribution under $\mathbb{P}_{\theta}$, and therefore  
$$\mathbb{P}_{\theta'}\left(\phi_{m,\alpha}\leq 0\right)\geq \delta-\varsigma.$$ 
Letting  $\varsigma$ go to $0$ enables to conclude.
 
\end{proof}

\begin{proof}[Proof of Proposition \ref{minoration_dependence}]

This lower bound for dependent gaussian covariates is proved through the same 
approach as Theorem \ref{minoration_arete}. We define the measure $\mu_{\rho}$ 
as in that proof. 
Under the hypothesis $H_0$, $Y$ is independent of $X$. We note $\Sigma$ the 
covariance matrix of $X$ and $\mathbb{E}_{0,\Sigma}$ stands for the 
distribution of $({\bf Y},{\bf X})$ under $H_0$ in order to emphasize the 
dependence on $\Sigma$. 
 
First, one has to upper bound the quantity $\mathbb{E}_{0,\Sigma}\left[L_{\mu_{\rho}}^2({\bf Y},{\bf X})\right]$. For the sake of simplicity, we make the hypothesis that every covariate $X_j$ 
has variance 1. If this is not the case, we only have to rescale these 
variables. The quantity $\corr(i,j)$ refers to the correlation between $X_i$ 
and $X_j$. As we only consider the case $k=1$, the set of models $m$ in 
$\mathcal{M}(1,p)$ is in correspondence with the set $\{1,\ldots , p\}$. 
 
\begin{eqnarray} 
\mathbb{E}_{0,\Sigma}\left(L_{\mu_{i,\zeta^1,\rho}}({\bf Y},{\bf X})L_{\mu_{j,\zeta^2,\rho}}({\bf Y},{\bf X})\right)=\left(\frac{\var(Y)}{\var(Y)-\corr(i,j)\lambda^2\zeta^1\zeta^2}\right)^n.\nonumber 
\end{eqnarray} 
 
When $i$ and $j$ are fixed, we upper bound the expectation of this quantity 
with respect to $\zeta^1$ and $\zeta^2$ by 
\begin{eqnarray} 
\mathbb{E}_{0,\Sigma}\left(L_{\mu_{i,\rho}}({\bf Y},{\bf X})L_{\mu_{j,\rho}}({\bf Y},{\bf X})\right) \leq 
\frac{1}{2} +\frac{1}{2}\left(\frac{\var(Y)}{\var(Y)-|\corr(i,j)|\lambda^2}\right)^n. \label{vraisemblance1_correle} 
\end{eqnarray} 
If $i\neq j$, $|\corr(i,j)|$ is smaller than $c$ and if $i=j$, $\corr(i,j)$ is exactly 
one. As a consequence, taking the expectation of 
(\ref{vraisemblance1_correle}) with respect to $i$ and $j$ yields the upper
bound 
\begin{equation} 
\mathbb{E}_{0,\Sigma}\left(L_{\mu_\rho}^2({\bf Y},{\bf X})\right) \leq \frac{1}{2} 
+\frac{1}{2}\left(\frac{1}{p}\left(\frac{\var(Y)}{\var(Y)-\lambda^2}\right)^n 
+ \frac{p-1}{p}\left(\frac{\var(Y)}{\var(Y)-c\lambda^2}\right)^n\right)\ .\label{vraisemblance2_correle} 
\end{equation}  
 
Recall that we want to constrain this quantity (\ref{vraisemblance2_correle}) 
to be smaller than $1+\eta^2$. In particular, this holds if 
the two following inequalities hold: 
\begin{eqnarray} 
\frac{1}{p}\left(\frac{\var(Y)}{\var(Y)-\lambda^2}\right)^n & \leq & 
\frac{1}{p} +\eta^2 \label{condition1_dependence}\\  
\frac{p-1}{p}\left(\frac{\var(Y)}{\var(Y)-c\lambda^2}\right)^n & \leq & 
\frac{p-1}{p}+ \eta^2\  .\label{condition2_dependence} 
\end{eqnarray} 
One then uses the inequality $\log(\frac{1}{1-x})\leq \frac{x}{1-x}$ which 
holds for any positive $x$ smaller than one. Condition (\ref{condition1_dependence}) holds if  
\begin{eqnarray} 
\frac{\lambda^2}{\var(Y)-\lambda^2} \leq \frac{1}{n}\log(1+p\eta^2)\ ,\label{condition3_dependence} 
\end{eqnarray} 
whereas Condition (\ref{condition2_dependence}) is implied by 
$$\frac{c\lambda^2}{\var(Y)-c\lambda^2} \leq \frac{1}{n}\log\left(1+\frac{p}{p-1}\eta^2\right)\ .$$ 
As $c$ is smaller than one and $\frac{p}{p-1}$ is larger than 1, this last inequality holds if 
\begin{eqnarray} 
\frac{\lambda^2}{\var(Y)-\lambda^2} \leq 
\frac{1}{nc}\log(1+\eta^2)\ . \label{condition4_dependence} 
\end{eqnarray} 
 
Gathering conditions (\ref{condition3_dependence}) and 
(\ref{condition4_dependence}) allows to conclude and to obtain the desired 
lower bound (\ref{minoration1_dependence}). 

\end{proof}
 
\begin{proof}[Proof of Proposition \ref{spatial}]
 
The sketch of the proof and the notations are analogous to the one in Proposition 
\ref{minoration_dependence}. The upper bound (\ref{vraisemblance1_correle}) 
still holds: 
\begin{eqnarray} 
\mathbb{E}_{0,\Sigma}\left(L_{\mu_{i,\rho}}({\bf Y},{\bf 
  X})L_{\mu_{j,\rho}}({\bf Y},{\bf X})\right) \leq 
\frac{1}{2} 
+\frac{1}{2}\left(\frac{\var(Y)}{\var(Y)-|\corr(i,j)|\lambda^2}\right)^n. \nonumber 
\end{eqnarray} 
Using the stationarity of the covariance function, we derive from 
(\ref{vraisemblance1_correle}) the following upper bound: 
\begin{eqnarray*} 
\mathbb{E}_{0,\Sigma}\left(L^2_{\mu_{\rho}}({\bf Y},{\bf X})\right) \leq 
\frac{1}{2} +\frac{1}{2p}\sum_{i=0}^{p-1}\left(\frac{\var(Y)}{\var(Y) - \lambda^2|\corr(0,i)|}\right)^n\ , 
\end{eqnarray*} 
where $\corr(0,i)$ equals $\corr(X_1,X_{i+1})$. As previously, we want to constrain this quantity to be smaller than $1+\eta^2$. In 
particular, this is implied if for any $i$ between $0$ and $p-1$: 
\begin{eqnarray*} 
\left(\frac{\var(Y)}{\var(Y)-\lambda^2|\corr(i,0)|}\right)^n\leq 1+ \frac{2p\eta^2|\corr(i,0)|}{\sum_{i=0}^{p-1}|\corr(i,0)|}. 
\end{eqnarray*} 
Using the inequality $\log(1+u)\leq u$, it is straightforward to show that this 
previous inequality holds if 
\begin{eqnarray*} 
\frac{\lambda^2}{\var(Y)-\lambda^2|\corr(i,0)|}\leq \frac{1}{n|\corr(i,0)|}\log\left(1+\frac{2p\eta^2|\corr(0,i)|}{\sum_{i=0}^{p-1}|\corr(i,0)|}\right). 
\end{eqnarray*} 
 
As $|\corr(i,0)|$ is smaller than one for any $i$ between $0$ 
and $p-1$, it follows that 
$\mathbb{E}_{0,\Sigma}\left(L^2_{\mu_{\rho}}({\bf Y},{\bf X})\right)$ is smaller than 
$1+\eta^2$ if 
\begin{eqnarray} 
\rho^2 & \leq & \bigwedge_{i=0}^{p-1} 
\frac{1}{n|\corr(i,0)|}\log\left(1+\frac{2p\eta^2|\corr(0,i)|}{\sum_{i=0}^{p-1}|\corr(i,0)|}\right). 
\nonumber 
\end{eqnarray} 
We now apply the convexity inequality $\log(1+ux)\geq u\log(1+x)$ which holds for any 
positive $x$ and any $u$ between $0$ and $1$ to obtain the condition  
\begin{eqnarray} 
\rho^2& \leq & \frac{1}{n}\log\left(1+\frac{2p\eta^2}{\sum_{i=0}^{p-1}|\corr(i,0)|}\right).\label{condition5_dependence} 
\end{eqnarray} 
 
It turns out we only have to upper bound the sum of $|\corr(i,0)|$ for the 
different types of correlation: 
\begin{enumerate} 
\item For $\corr(i,j)= \exp(-w|i-j|_p)$, the sum is clearly bounded by 
$1+2\frac{e^{-w}}{1-e^{-w}}$ and 
  Condition (\ref{condition5_dependence}) simplifies as 
$$\rho^2\leq 
\frac{1}{n}\log\left(1+2p\eta^2\frac{1-e^{-w}}{1+e^{-w}}\right). 
 $$ 
\item if $\corr(i,j) = \left(1+|i-j|_p\right)^{-t}$ for $t$ strictly larger 
  than one, then $\sum_{i=0}^{p-1}|\corr(i,0)|\leq 1 +\frac{2}{t-1}$ and 
  Condition (\ref{condition5_dependence}) simplifies as 
$$\rho^2\leq \frac{1}{n}\log\left(1+\frac{2p(t-1)\eta^2}{t+1}\right).$$ 
\item if $\corr(i,j) = \left(1+|i-j|_p\right)^{-1}$  then $\sum_{i=0}^{p-1}|\corr(i,0)|\leq 1 +2\log(p-1)$ and 
  Condition (\ref{condition5_dependence}) simplifies as 
$$\rho^2\leq \frac{1}{n}\log\left(1+\frac{2p\eta^2}{1+2\log(p-1)}\right).$$ 
\item if $\corr(i,j) = \left(1+|i-j|_p\right)^{-t}$ for $0<t<1$, then $$\sum_{i=0}^{p-1}|\corr(i,0)|\leq 1 +\frac{2}{1-t}\left[\left(\frac{p}{2}\right)^{1-t}-1\right]\leq \frac{2}{1-t}\left(\frac{p}{2}\right)^{1-t}$$ and 
  Condition (\ref{condition5_dependence}) simplifies as 
  $$\rho^2\leq \frac{1}{n}\log\left(1+p^t2^{1-t}(1-t)\eta^2\right).$$ 
\end{enumerate} 
 \end{proof}
 
\begin{proof}[Proof of Proposition \ref{vitesse_minimax_ellipsoide}]
 
For each dimension $D$ between $1$ and $p$, we define $r_D^2=\rho^2_{D,n}\wedge 
a_D^2R^2$. Let us fix some $D\in \{1,\ldots,p\}$. Since $r_D^2\leq a_D^2$ and since the $a_j$'s are 
non increasing, $$\sum_{j=1}^D\frac{\var(Y|X_{m_{j-1}}) -\var(Y|X_{m_{j}})} 
{a_j^2}\leq \var(Y|X)R^2,$$ for all $\theta\in S_{m_D}$ such that 
$\frac{\|\theta\|^2}{\var(Y)-\|\theta\|^2} = r^2_D$. Indeed, $\|\theta\|^2 = 
\sum_{j=1}^D\var(Y|X_{m_{j-1}}) -\var(Y|X_{m_{j}})$ and $\var(Y)-\|\theta\|^2= 
\var(Y|X)$. As a consequence,  
\begin{eqnarray*} 
\left\{\theta\in S_{m_D}, \frac{\|\theta\|^2}{\var(Y)-\|\theta\|^2}=r^2_D\right\}\subset \left\{\theta\in \mathcal{E}_a(R), \frac{\|\theta\|^2}{\var(Y)-\|\theta\|^2}\geq r^2_D\right\}\ . 
\end{eqnarray*} 
Since $r_D\leq \rho_{D,n}$, we deduce from Proposition \ref{minoration_1test} 
that 
\begin{eqnarray*} 
\beta_{\Sigma}\left(\left\{\theta\in \mathcal{E}_a(R), \frac{\|\theta\|^2}{\var(Y)-\|\theta\|^2}\geq r_D^2 \right\}\right)\geq \delta\ .
\end{eqnarray*} 
The first result of Proposition \ref{vitesse_minimax_ellipsoide} follows by 
gathering these lower bounds for all $D$ between $1$ and $p$. 
 
Moreover, $\rho^2_{i,n} $ is defined in Proposition \ref{minoration_1test} as 
$\rho^2_{i,n}=\sqrt{2}\left[\sqrt{\mathcal{L}(\eta)}\wedge\frac{\mathcal{L}(\eta)}{\sqrt{\log 
      2}}\right]\frac{\sqrt{i}}{n}$. If $\alpha+\delta\leq 47\%$, it is 
straightforward to show that $\rho^2_{i,n}\geq \frac{\sqrt{i}}{n}$. 
 
 \end{proof}

\begin{proof}[Proof of Proposition \ref{minimax_lower_nested_linear}]
 
We first need the following Lemma. 
\begin{lemma}\label{minimax_simultaneous_detection} 
We consider $(I_j)_{j\in \mathcal{J}}$ a partition of $\mathcal{I}$. For each 
$j\in \mathcal{J}$ let $p(j)=|I_{j}|$. For any $j\in \mathcal{J}$, we define $\Theta_j$ as the set of 
$\theta\in \mathbb{R}^{p}$ such that their support is included in $I_j$. 
For any sequence of positive weights $k_j$ such that 
\begin{eqnarray*} 
\sum_{j\in \mathcal{J}}k_j=1, 
\end{eqnarray*} 
it holds that  
\begin{eqnarray*} 
\beta_I\left(\bigcup_{j\in \mathcal{J}}\left\{\theta\in \Theta_j,\frac{\|\theta\|^2}{\text{\emph{var}}(Y)-\|\theta\|^2}=r_j^2   \right\}\right)\geq \delta \ , 
\end{eqnarray*} 
if for all $j\in \mathcal{J}$, $r_j\leq \rho_{p(j),n}(\eta/\sqrt{k_j})$, where the function
$\rho_{p(j),n}$ is defined by (\ref{rho_1test}). 
\end{lemma} 
 For all $j\geq 0$ such that 
$2^{j+1}-1\in \mathcal{I}$ (i.e. for all $j\leq J$ where $J = 
\log(p+1)/\log(2)-1\ $), let $\bar{S}_j$ be the linear span of the $e_k$'s for 
$k\in \{2^j,\ldots,2^{j+1}-1\}$. Then, dim$(\bar{S}_j)=2^j$ and 
$\bar{S}_j\subset S_{m_D}$ for $D=D(j)=2^{j+1}-1$. It is straightforward to show that 
\begin{eqnarray*} 
\bigcup_{j=0}^J \bar{S}_j[r_{D(j)}] \subset \bigcup_{j=0}^J 
S_{m_{D(j)}}[r_{D(j)}]\subset \bigcup_{D=1}^p S_{m_D}[r_D]\ , 
\end{eqnarray*} 
where $\bar{S}_j[r_D(j)] := \left\{\theta\in \bar{S}_j 
,\frac{\|\theta\|^2}{\var(Y)-\|\theta\|^2}=r_{D(j)}^2 
\right\}$ and $ S_{m_D}[r_D] :=\left\{\theta\in S_{m_D} 
,\frac{\|\theta\|^2}{\var(Y)-\|\theta\|^2}=r_{D}^2 
\right\}$.  
 
Let choose $\mathcal{J}=\left\{1,\ldots,J\right\}$. For any $j\in\mathcal{J}$, we define $I_j=\left\{2^j,2^j+1,\ldots,2^{j+1}-1\right\}$.
Applying Lemma \ref{minimax_simultaneous_detection} with $k_j := 
\left[(j+1)R(p)\right]^{-1}$ where $R(p):=\sum_{k=0}^J 1/(k+1) $ we get 
\begin{eqnarray*} 
\beta_{I}\left(\bigcup_{D=1}^p \left\{\theta\in S_{m_D}, \frac{\|\theta\|^2}{\var(Y)-\|\theta\|^2}=r_D^2 \right\}\right)\geq \delta\ , 
\end{eqnarray*} 
if for all those $D=D(j)$ 
\begin{eqnarray*} 
r_D^2\leq \sqrt{\log(1+2\eta^2/k_j)}\left(1\wedge \frac{\sqrt{\log(1+2\eta^2/k_j)}}{\sqrt{2\log 2}}\right)\frac{\sqrt{D}}{n}. 
\end{eqnarray*} 
For $D=D(j)$, this last quantity is lower bounded by
\begin{eqnarray}\label{condition_nested_linear} 
 \sqrt{\log(1+2\eta^2/k_j)}\left(1\wedge 
 \frac{\sqrt{\log(1+2\eta^2/k_j)}}{\sqrt{2\log 2}}\right)\frac{\sqrt{D}}{n} \geq\hspace{4cm} \\ \hspace{4cm} \sqrt{\log(1+2\eta^2(j+1)R(p))}\left(1\wedge 
 \frac{\sqrt{\log(1+2\eta^2)}}{\sqrt{2\log 2}}\right) \frac{2^{j/2}}{n}\nonumber\ . 
\end{eqnarray} 
It remains to check that (\ref{condition_nested_linear}) is larger than $\bar{\rho}_{D(j),n}$. Using 
$j+1 = \log(D+1)/\log(2)\geq \log(D+1)$, we get 
$2^{j/2}\geq \sqrt{D/2}$. Thanks to the convexity inequality 
$\log(1+ux)\geq u\log(1+x)$, 
which holds for any $x>0$ and any $u\in]0,1]$, we obtain 
\begin{eqnarray*} 
\sqrt{\log(1+2\eta^2(j+1)R(p))}2^{j/2}& \geq & 
\sqrt{D/2}\left(\eta\sqrt{2R(p)}\wedge 1\right)\sqrt{\log\left[1+\log(D+1)\right]}\\ 
&\geq & \left( (\eta\sqrt{2})\wedge 1 \right)\sqrt{\log \log (D+1)}\sqrt{D/2},\\ 
& \geq & \frac{1}{\sqrt{2}}\left( 1\wedge \sqrt{\log(1+2\eta^2)}\right)\sqrt{\log \log (D+1)}\sqrt{D}\ , 
\end{eqnarray*} 
as $R(p)$ is larger than one for any $p\geq 1$. All in all, we get the lower 
bound 
\begin{eqnarray*} 
& \sqrt{\log(1+2\eta^2(j+1)^2R(p))} & \left(1\wedge 
 \frac{\sqrt{\log(1+2\eta^2)}}{\sqrt{2\log 2}}\right) \frac{2^{j/2}}{n} \\
& &\geq \frac{1}{2\sqrt{\log(2)}}\left(1\wedge 
 \log(1+2\eta^2)\right)\sqrt{\log\log (D+1)}\frac{\sqrt{D}}{n} = \bar{\rho}_{D,n}^2\ . 
\end{eqnarray*} 
Thus, if for all $1\leq D\leq p$, $r_D^2$ is smaller than $\bar{\rho}^2_{D,n}$, it holds that 
\begin{eqnarray*} 
\beta_{I}\left(\bigcup_{D=1}^p \left\{\theta\in S_{m_D}, 
\frac{\|\theta\|^2}{\var(Y)-\|\theta\|^2}=r_D^2  \right\}\right)\geq \delta\ . 
\end{eqnarray*} 
\end{proof} 
 
\begin{proof}[Proof of Lemma \ref{minimax_simultaneous_detection}]
 
Using a similar approach to the proof of Theorem \ref{minoration_arete}, we 
know that for each $r_j\leq \tilde{\rho}_j(\eta/\sqrt{k_j})$ there exists some 
measure $\mu_j$ over 
$$\Theta_j[r_j] := \left\{\theta\in 
\Theta_j,\frac{\|\theta\|^2}{\var(Y)-\|\theta\|^2} =r_j^2 \right\}$$ 
such that 
\begin{eqnarray}\label{expression1_detectionsimutanee} 
\mathbb{E}_0\left[L^2_{\mu_j}(Y,X)\right]\leq 1+\eta^2/k_j\ . 
\end{eqnarray} 
We now define a probability measure $\mu=\sum_{j\in \mathcal{J}}k_j\mu_j$ over 
$\bigcup_{j\in\mathcal{J}}\Theta_j[r_j]$. $L_{\mu_j}$ refers to the density of 
$\mathbb{P}_{\mu_j}$ with respect to $\mathbb{P}_0$. Thus, 
\begin{eqnarray*} 
L_{\mu}(Y) = \frac{\text{d}\mathbb{P}_{\mu}}{\text{d}\mathbb{P}_0}({\bf Y},{\bf X}) = \sum_{j\in 
  \mathcal{J}}k_jL_{\mu_j}({\bf Y},{\bf X})\ , 
\end{eqnarray*} 
and 
\begin{eqnarray*} 
\mathbb{E}_0\left[L^2_{\mu}({\bf Y},{\bf X})\right] = \sum_{j,j'\in 
  \mathcal{J}}k_jk_{j'}\mathbb{E}_0\left[L_{\mu_j}({\bf Y},{\bf X}) L_{\mu_{j'}}({\bf Y},{\bf X}) \right]\ . 
\end{eqnarray*} 
Using expression (\ref{expression_vraisemblance}), it is straightforward 
to show that if $j\neq j'$, then $$\mathbb{E}_0\left[L_{\mu_j}({\bf Y},{\bf X}) L_{\mu_{j'}}({\bf Y},{\bf X}) \right]=1.$$ This follows from the fact that 
the sets $\Theta_j$ and $\Theta_{j'}$ are orthogonal with respect to the inner product (\ref{definition_produit}). Thus, 
\begin{eqnarray*} 
\mathbb{E}_0\left[L_{\mu}({\bf Y},{\bf X})\right] = 1 + 
\sum_{j\in\mathcal{J}}k_j^2\left(\mathbb{E}_0\left[L^2_{\mu_j}({\bf Y},{\bf X})\right]-1 \right)\leq 1+\eta^2 
\end{eqnarray*} 
thanks to (\ref{expression1_detectionsimutanee}). Using the 
argument (\ref{argument_principal_minimax}) as in the proof of Theorem 
\ref{minoration_arete} allows to conclude. 
 \end{proof}

\begin{proof}[Proof of Proposition \ref{minoration_minimax_collection_ellipsoide}]
 
First of all, we only have to consider the case where the covariance matrix of 
$X$ is the identity. If this is not the case, one only has to apply 
Gram-Schmidt process to $X$ and thus obtain a vector $X'$ and a new basis for 
$\mathbb{R}^p$ which is orthonormal. We refer to the beginning of Section 
\ref{ellipsoides} for more details. 
 
Like the previous bounds for ellipsoids, we adapt the approach of Section 6 in 
Baraud \cite{baraudminimax}. We use the same notations as in proof of Proposition 
\ref{vitesse_minimax_ellipsoide}. Let $D^*(R) \in \{1,\ldots, p\}$ an integer 
which achieves the supremum of $\bar{\rho}_D^2\wedge (R^2a_D^2)=\bar{r}_D^2$.  
As in proof of Proposition \ref{vitesse_minimax_ellipsoide}, 
for any $R>0$, 
\begin{eqnarray*} 
\left\{\theta\in S_{m_{D^*(R)}}, \frac{\|\theta\|^2}{\var(Y)-\|\theta\|^2}=r^2_{D^*(R)} \right\}\subset \left\{\theta\in 
\mathcal{E}_a(R), \frac{\|\theta\|^2}{\var(Y)-\|\theta\|^2}\geq r^2_{D^*(R)} \right\}\ . 
\end{eqnarray*} 
When $R$ varies, $D^*(R)$ describes $\{1,\ldots, p\}$. Thus, we obtain  
\begin{eqnarray*} 
\bigcup_{1\leq D\leq p}\left\{\theta\in 
S_{m_D}, \frac{\|\theta\|^2}{\var(Y)-\|\theta\|^2}=r_D^2 \right\}&= &
\bigcup_{R>0}\left\{\theta\in
S_{m_{D^*(R)}}, \frac{\|\theta\|^2}{\var(Y)-\|\theta\|^2}=
r^2_{D^*(R)}\right\} \\ 
& \subset & \bigcup_{R>0}\left\{\theta\in \mathcal{E}_a(R)
\frac{\|\theta\|^2}{\var(Y)-\|\theta\|^2}\geq r^2_{D^*(R)} \right\}\ , 
\end{eqnarray*}
and the result follows from proposition \ref{minimax_lower_nested_linear}.
\end{proof}

\appendix
\section*{Appendix}
 \begin{proof}[Proof of Proposition \ref{niveau}]
The test associated with Procedure $P_1$ corresponds to a Bonferroni procedure. Hence, we prove that its size is less than $\alpha$ by arguing as follows: 
let $\theta$ be an element of $S_V$ (defined in Section \ref{presentation}),  
\begin{eqnarray*} 
\mathbb{P}_{\theta} (T_{\alpha} > 0) 
& \leq & \sum_{m \in \mathcal{M}}  \mathbb{P}_{\theta}  \left(\phi_m({\bf 
  Y},{\bf X}) - 
\bar{F}^{-1}_{D_m,N_m}(\alpha_m) > 0 \right), 
\end{eqnarray*} 
where $\phi_m({\bf Y},{\bf X})$ is defined in (\ref{definition_phi}).
The test is rejected if for some model $m$, $\phi_{m}({\bf Y},{\bf X})$ is larger than $\bar{F}^{-1}_{D_m,N_m}(\alpha_m)$. 
As $\theta$ belongs to $S_V$, $\Pi_{V\cup m}{\bf Y} - \Pi_V{\bf Y}=\Pi_{V\cup 
      m}\boldsymbol{\epsilon} -\Pi_V \boldsymbol{\epsilon}$ 
      and ${\bf Y} - \Pi_{V\cup m}{\bf Y} = \boldsymbol{\epsilon} - \Pi_{V\cup 
	m}\boldsymbol{\epsilon}$. Then, the quantity $\phi_m({\bf 
  Y},{\bf X})$ is equal to  
$$\phi_m({\bf 
  Y},{\bf X}) =  \frac{N_m\|\Pi_{V\cup m}{\boldsymbol{\epsilon}} - 
    \Pi_V{\boldsymbol{\epsilon}}\|^2_n}{D_m\|{ \boldsymbol{\epsilon}} - \Pi_{V\cup 
    m}{\boldsymbol{\epsilon}}\|^2_n}.$$ 
Because $\boldsymbol{\epsilon}$ is independent of ${\bf X}$, the distribution of 
  $\phi_m({\bf Y},{\bf X})$ conditionally to ${\bf X}$ is a Fisher
  distribution with $D_m$ and 
  $N_m$ degrees of freedom. As a consequence, $\phi_{m,\alpha_m}({\bf Y},{\bf X})$ 
  is a Fisher test with $D_m$ and $N_m$ degrees of freedom. It 
  follows that: 
$$\mathbb{P}_{\theta}(T_{\alpha} > 0) 
 \leq \sum_{m \in \mathcal{M}}\alpha_m\leq \alpha. $$ 
\vspace{0.5cm}

The test associated with Procedure $P_2$  has the property to 
be of size exactly $\alpha$. More precisely, for any $\theta\in S_V$, we have that   
$$\mathbb{P}_{\theta}(T_{\alpha} >0|{\bf X}) = \alpha\, \, \, \, \, \, {\bf 
  X}\ \text{a.s.}\ .$$ 
The result follows from the fact that  $q_{{\bf X},\alpha}$ satisfies 
$$\mathbb{P}_{\theta}\left(\left. \sup_{m \in \mathcal{M}} \left\{  \frac{N_m \| 
      \Pi_{V\cup m} (\boldsymbol{\epsilon}) -  
  \Pi_{V}(\boldsymbol{\epsilon})  \|_n^2 }{ D_m\|\boldsymbol{\epsilon} -  
  \Pi_{V\cup m}(\boldsymbol{\epsilon}) \|_n^2 }  - 
\bar{F}^{-1}_{D_m, N_m}\left(q_{{\bf X},\alpha}\right)\right\} >0 \right|{\bf X}    \right)= \alpha,$$ 
and that for any $\theta\in S_V$, $\Pi_{V\cup m}{\bf Y} - \Pi_V{\bf Y}=\Pi_{V\cup 
      m}\boldsymbol{\epsilon} -\Pi_V \boldsymbol{\epsilon}$ 
      and ${\bf Y} - \Pi_{V\cup m}{\bf Y} = \boldsymbol{\epsilon} - \Pi_{V\cup m}\boldsymbol{\epsilon}$. 
  
 \end{proof}

\begin{proof}[Proof of Proposition \ref{puissance_comparaison}]
Let come back to the  
definitions of $T^1_\alpha$ and $T^2_{\alpha}$: 
\begin{eqnarray} 
T^1_{\alpha}({\bf X}, {\bf Y}) & = & \sup_{m\in \mathcal{M}}\left\{\phi_m({\bf Y},{\bf X})-\bar{F}^{-1}_{D_m,N_m}(\alpha/|\mathcal{M}|)\right\}\nonumber\\ 
T^2_{\alpha}({\bf X},{\bf Y}) & = & \sup_{m\in \mathcal{M}}\left\{\phi_m({\bf Y},{\bf 
  X})-\bar{F}^{-1}_{D_m,N_m}(q_{\bf X,\alpha})\right\}\nonumber 
\end{eqnarray} 
Conditionally on ${\bf X}$, the size of $T^1_\alpha$ is smaller than $\alpha$, 
whereas  the size $T^2_{\alpha}$ is exactly $\alpha$. As a consequence $q_{\bf X,\alpha}\geq 
\alpha/|\mathcal{M}|$ as the statistics $T^1_{\alpha}$ and $T^2_{\alpha}$ differ only through these 
quantities. Thus, $T^2_{\alpha}({\bf X},{\bf Y})\geq 
T^1_{\alpha}({\bf X}, {\bf Y})$, $({\bf X},{\bf Y})$ almost surely and the 
result (\ref{comparaison_puissance}) follows.
\end{proof}

\section*{Acknowledgements}
We gratefully thank Sylvie Huet and Pascal Massart for many fruitful discussions.

\addcontentsline{toc}{section}{References}
\bibliographystyle{elsart-harv}
\bibliography{tests}

\end{document}